\documentclass[a4paper,11pt,twoside,reqno]{amsart}

\usepackage[utf8]{inputenc}
\usepackage[plainpages=false,pdfpagelabels=true]{hyperref}
\usepackage{amssymb,amsthm}
\usepackage[margin=1in]{geometry}
\usepackage{slashed}
\usepackage{color}
\usepackage[foot]{amsaddr}
\newtheorem{Satz}{Theorem}[section]
\newtheorem{Prop}[Satz]{Proposition}
\newtheorem{Lem}[Satz]{Lemma}
\newtheorem{Thm}[Satz]{Theorem}
\newtheorem{Cor}[Satz]{Corollary}
\theoremstyle{definition}
\newtheorem{Dfn}[Satz]{Definition}
\newtheorem{Bem}[Satz]{Remark}

\allowdisplaybreaks[1]

\renewcommand{\epsilon}{\varepsilon}

\newcommand{\R}{\ensuremath{\mathbb{R}}}

\numberwithin{equation}{section}

\newcommand{\p}{\slashed{\partial}}

\title{Eigenvalue estimates on weighted manifolds}
\author{Volker Branding} 
\address{University of Vienna, Faculty of Mathematics
Oskar-Morgenstern-Platz 1, 1090 Vienna, Austria}
\email{volker.branding@univie.ac.at}

\author{Georges Habib}
\address{Lebanese University, Faculty of Sciences {\bf II}, Department of Mathematics, 
P.O. Box 90656 Fanar-Matn, Lebanon and Universit\'e de Lorraine, CNRS, IECL, F-54506, Nancy, France}
\email{ghabib@ul.edu.lb}

\date{\today}

\thanks{
The first author gratefully acknowledges the support of the Austrian Science Fund (FWF) through the START-project (Y963) of Michael Eichmair and the project ``Geometric Analysis of Biwave Maps`` (P34853).
The second author gratefully acknowledges the support of the Humboldt Foundation and also 
wishes to thank Sylvie Paycha, the University of Potsdam and the MPI Bonn.
}
\begin{document}

\begin{abstract}
We derive various eigenvalue estimates for the Hodge Laplacian acting on differential forms on weighted Riemannian manifolds.
Our estimates unify and extend various results from the literature and
 provide a number of geometric applications. In particular, we derive an inequality
which relates the eigenvalues of the Jacobi operator for \(f\)-minimal hypersurfaces
and the spectrum of the Hodge Laplacian.
\end{abstract} 

\maketitle

\section{Introduction and results}

The aim of spectral geometry is to obtain deep insights into both topology and geometry
of Riemannian manifolds from the spectrum of differential operators.
However, in the case of an arbitrary Riemannian manifold, it is in general not possible to explicitly
calculate the spectrum of a given differential operator.
For this reason, it is very important to obtain characterizations of the eigenvalues of such operators in terms of geometric data as for example curvature.

The most prominent spectral problem one can study is the case of the Laplace operator acting on functions
for which we recall several important results.
Let \(\Omega\subset M\) be a bounded domain of a Riemannian manifold $M$  with smooth boundary \(\partial\Omega\) 
and let \(\lambda\) be an eigenvalue of the Dirichlet problem
\begin{equation}\label{eq:euclideandirichlet}
\left\{\begin{array}{ll}
\Delta u=\lambda u\,\,  &\text{on}\,\, \Omega\\
u=0 \,\,&\text{on}\,\, \partial\Omega
\end{array}\right.
\end{equation}
where \(u\in C^\infty(M)\) and \(\Delta\) is the Laplace-Beltrami operator.
From general spectral theory, we know that the spectrum of the Dirichlet problem \eqref{eq:euclideandirichlet} is discrete and that the eigenvalues satisfy
\begin{align*}
0<\lambda_1\leq\lambda_2\leq\lambda_3\leq\ldots\leq\lambda_j\leq\ldots\to\infty. 
\end{align*}

In the case that the manifold $M$ is the Euclidean space with the flat metric \((\mathbb{R}^n,{\rm can})\),
Payne, P\'olya and Weinberger \cite{PGW1,PGW2} derived the following
so-called \emph{universal inequalities} for the Dirichlet problem of the Laplace operator \eqref{eq:euclideandirichlet}.
Namely, they show that for any integer $k\geq 1$ 
\begin{align}
\label{eq:inequality-intro}
\lambda_{k+1}-\lambda_k\leq\frac{4}{nk}\sum_{i=1}^k\lambda_i.
\end{align}
An immediate consequence of \eqref{eq:inequality-intro} is the following inequality,
$$\lambda_{k+1}\leq (1+\frac{4}{n})\lambda_k,$$
which allows to get an estimate of all the eigenvalues if an upper bound is known on the lowest one $\lambda_1$. Another kind of such inequality has been further developed by Hile and Protter (see \cite{HP}) and by Yang in \cite{Y} who 
established the following estimate
\begin{align}
\label{eq:recursionintroduction}
\sum_{i=1}^k(\lambda_{k+1}-\lambda_i)^2\leq\frac{4}{n}\sum_{i=1}^k(\lambda_{k+1}-\lambda_i)\lambda_i.
\end{align}
This inequality improves all the previous estimates. In \cite{CY}, Cheng and Yang were able to turn the recursion formula \eqref{eq:recursionintroduction}
into an estimate for the \((k+1)\)-th eigenvalue of \eqref{eq:euclideandirichlet} in terms of some power of $k$  as follows 
\begin{align}
\lambda_{k+1}\leq C_0(n,k)k^\frac{2}{n}\lambda_1,
\end{align}
where \(C_0(n,k)\) is a constant only depending on \(n\) and \(k\).
For more details on the spectral properties of the Dirichlet problem \eqref{eq:euclideandirichlet}
we refer to the lecture notes \cite{Ash}.

Many of these formulas have been generalized to the case of a Riemannian manifold. In this setup, Chen and Cheng \cite{CC} obtained an extrinsic estimate for eigenvalues of the Dirichlet problem \eqref{eq:euclideandirichlet} of the Laplacian on a complete Riemannian manifold \((M^n,g)\) that is isometrically 
immersed in an \((n+m)\)-dimensional Euclidean space \(\R^{n+m}\) (see also \cite{EHI}). The estimate is 
\begin{align}\label{eq:estimatemean}
\sum_{i=1}^k(\lambda_{k+1}-\lambda_i)^2\leq
\frac{4}{n}\sum_{i=1}^k(\lambda_{k+1}-\lambda_i)(\lambda_i+\frac{n^2}{4}\sup_\Omega |H|^2).
\end{align}

Here, $H:=\frac{1}{n}{\rm trace}({\bf II})$ denotes the mean curvature of the immersion and ${\bf II}$ is the second fundamental form. As a consequence, they deduce an upper bound for the $(k+1)$-th eigenvalue in terms of a power of $k$ that involves the mean curvature. 

In this article, we focus on the study of geometric differential operators
on an interesting class of Riemannian manifolds, the so-called \emph{Bakry-Emery manifolds}.
A Bakry-Emery manifold is a triple \((M,g,d\mu_f)\), where instead of the Riemannian measure \(dv_g\)
one considers the weighted measure
\(d\mu_f:=e^{-f}dv_g\) with \(f\in C^\infty(M)\). One often refers to this kind of manifolds as 
\emph{weighted manifolds}. While the initial motivation to study such kind of manifolds was to model diffusion processes \cite{BE}
they have by now become famous in the study of self-similar solutions of the Ricci flow,
the so-called \emph{Ricci solitons}. Also they appear in the analysis of \emph{shrinkers},
which represent a special class of solutions of the mean curvature flow. We refer to \cite[Section 2]{CM} and \cite{IRS} for more details. 

Due to the presence of the weight, differential operators on weighted manifolds
contain additional contributions and also their spectrum is different compared
to the case of standard Riemannian manifolds.
For example, the Laplace operator acting on functions on a weighted manifold 
is defined as follows
\begin{align}
\label{eq:defdriftlaplace}
\Delta_f:=\Delta+\nabla_{df},
\end{align}
where \(\Delta=\delta d\) is the standard Laplace-Beltrami operator.
 A direct computation shows that \(\Delta_f\), usually called {\it drift Laplacian}, is elliptic and self-adjoint (if $M$ is compact)  with respect to the weighted measure \(e^{-f}dv_g\). Therefore, its spectrum is discrete and consists of an increasing sequence of eigenvalues  $$0<\lambda_{1,f}\leq\lambda_{2,f}\leq \dots\leq \lambda_{j,f}\leq\dots\to \infty.$$
As in the standard case, the eigenvalue $0$ corresponds to constant functions. Several results have been obtained on the spectrum of this operator (see for example \cite{MD, ML, ML1})  
of which we give a non-exhaustive list below. In particular, when a compact Bakry-Emery manifold $(M^n,g)$ is isometrically immersed into the Euclidean space,
Xia and Xu \cite{XX} established the following recursion formula 
\begin{align*}
\sum_{i=1}^k(\lambda_{k+1,f}-\lambda_{i,f})^2\leq
\frac{1}{n}\sum_{i=1}^k(\lambda_{k+1,f}-\lambda_{i,f})
\big(4\lambda_{i,f}
+4\sup_\Omega|\nabla f|\sqrt{\lambda_{i,f}}+n^2\sup_{\Omega}|H|^2
+\sup_\Omega|\nabla f|^2)\big)
\end{align*}
which clearly reduces to Inequality \eqref{eq:estimatemean} when $f$ is zero. Also  several explicit calculations for the drift Laplacian on self shrinkers were carried out in \cite{CP, Z1, BK}. 

This article is devoted to the study of the Hodge Laplacian acting on differential forms on weighted manifolds, which we will often refer to as \emph{drifting Hodge Laplacian}, and also denote it by \(\Delta_f\) (see Subsection \ref{sec:laplacweighted} for the definition). As for functions, the drifting Hodge Laplacian on differential $p$-forms has a discrete spectrum that entirely consists of nonnegative eigenvalues denoted by $(\lambda_{i,p,f})_i$.
In particular, we will derive various eigenvalue estimates as well as different recursion formulas which will characterize the corresponding spectrum.

In the following, we present the main results of this article. 
First, we extend the results of \cite{Asa,S} to the case of weighted manifolds:
\begin{Thm} \label{thm:inequality} 
Let $(M^n,g, d\mu_f=e^{-f}dv_g)$ be a compact Bakry-Emery manifold that is isometrically immersed into the Euclidean space $\R^{n+m}$. For $p=1,\ldots,n$, we let $\lambda_{1,p,f}$ be  the first nonnegative eigenvalue of $\Delta_f$ on $p$-forms and $\lambda'_{1,p,f}$   the first positive eigenvalue of $\Delta_f$ on exact $p$-forms. Then, we have the following estimate 
\begin{align}
\lambda_{1,p,f}-\lambda_{1,p-1,f}\geq \frac{1}{p}\mathop {\rm inf}\limits_M\left(\mathfrak{B}^{[p]}+T_f^{[p]}-\sum_{s=1}^m({\bf II}^{[p]}_{\nu_s})^2\right).
\end{align}
Also, we have the inequality
\begin{align}
\lambda'_{1,p,f}\leq \frac{1}{{\rm Vol}_f(M)}\int_M \left(pn|H|^2-\frac{p(p-1)}{n(n-1)}{\rm Scal}^M+\frac{p}{n}|df|^2\right) d\mu_f,
\end{align}
where ${\rm Scal}^M$ denotes the scalar curvature of 
\(M\).
\end{Thm}

The tensor $\mathfrak{B}^{[p]}+T_f^{[p]}$ appearing in the first statement of Theorem \ref{thm:inequality}
is the so-called $p$-Ricci tensor, see \eqref{eq:Bochner} for the precise definition.
Also, ${\bf II}^{[p]}_{\nu_s}$ is some canonical extension of the second fundamental form ${\bf II}$ of the immersion to differential $p$-forms in the direction of the normal vector field \(\nu_s\), 
see \eqref{dfn:sff-forms} for more details.

\begin{Bem} 
\begin{enumerate} 
\item It is well-known that $\lambda_{1,f}=\lambda'_{1,1,f}$ where $\lambda_{1,f}$ is the first positive eigenvalue of $\Delta_f$ on functions. Hence we get the estimate for $\lambda_{1,f}$ \begin{equation}\label{eq:equalitycase}
\lambda_{1,f}\leq \frac{1}{{\rm Vol}_f(M)}\int_M \left(n|H|^2+\frac{1}{n}|df|^2\right) d\mu_f,
\end{equation}
which is the same estimate as in \cite{BCP}. 
\item It was shown in \cite{JMZ} that when $M$ is an embedded shrinker of revolution in $\R^3$ such that the intersection of $M$ with some sphere has only two connected components and such that $M$ is symmetric with respect to reflection across
the axis of revolution, then equality is attained in \eqref{eq:equalitycase}. Therefore, the inequality in Theorem \ref{thm:inequality} is attained for $p=1$.
\end{enumerate}
\end{Bem}

In addition to the eigenvalue estimates presented in Theorem \ref{thm:inequality}, we 
establish the following recursion formulas for the eigenvalues of the drifting Hodge Laplacian
on weighted manifolds:

\begin{Thm} 
\label{thm:recursionformulas} 
Let $X:(M^n,g)\to (\mathbb{R}^{n+m},{\rm can})$ be an isometric immersion. 
For any $p\in\{0,\ldots,n\}$, the eigenvalues of the drifting Hodge Laplacian $\Delta_f$ acting on $p$-forms on a domain $\Omega$ of $M$ with Dirichlet boundary conditions satisfy for any $k\geq 1$
\begin{eqnarray*}
\sum_{i=1}^k(\lambda_{k+1,p,f}-\lambda_{i,p,f})^\alpha&\leq& \frac{4}{n}\sum_{i=1}^k(\lambda_{k+1,p,f}-\lambda_{i,f})^{\alpha-1}\big(\lambda_{i,p,f}
-\int_\Omega\langle(\mathfrak{B}^{[p]}+T_f^{[p]})\omega_i,\omega_i\rangle d\mu_f\\&&
+\frac{n^2}{4}\int_\Omega|H|^2|\omega_i|^2d\mu_f
-\frac{1}{4}\int_\Omega (2\Delta f+|df|^2)|\omega_i|^2d\mu_f\big),
\end{eqnarray*}
for $\alpha\leq 2$. Also, we have 
\begin{eqnarray*}
\sum_{i=1}^k(\lambda_{k+1,p,f}-\lambda_{i,p,f})^\alpha&\leq& \frac{2\alpha}{n}\sum_{i=1}^k(\lambda_{k+1,p,f}-\lambda_{i,p,f})^{\alpha-1}\big(\lambda_{i,p,f}-\int_\Omega\langle(\mathfrak{B}^{[p]}+T_f^{[p]})\omega_i,\omega_i\rangle d\mu_f\\&&
+\frac{n^2}{4}\int_\Omega|H|^2|\omega_i|^2d\mu_f
-\frac{1}{4}\int_\Omega (2\Delta f+|df|^2)|\omega_i|^2d\mu_f\big),
\end{eqnarray*}
for $\alpha\geq 2$. 
Here, \(\omega_i\) represents the \(i\)-th eigenform of \(\Delta_f\).
\end{Thm}

\begin{Bem}
In Theorem \ref{thm:inequality} and Theorem \ref{thm:recursionformulas} the tensor \(\mathfrak{B}^{[p]}\)
appears in the estimates. It is possible to express this tensor in terms of extrinsic geometric data
via equation \eqref{eq:bp} but as this would further complicate the presentation of the eigenvalue estimates
we do not make use of this option.
\end{Bem}

In order to eliminate the dependence on $\omega_i$ in the statements of Theorem \ref{thm:recursionformulas}, we take the infimum over $\Omega$ of the smallest eigenvalue of the endomorphsim 
$\mathfrak{B}^{[p]}+T_f^{[p]}$ which we denote by $\delta_1$ to 
deduce that 
\begin{Cor} 
\label{cor:recursionformulas}
Let $X:(M^n,g)\to (\mathbb{R}^{n+m}, {\rm can})$ be an isometric immersion. For any $p\in\{0,\ldots,n\}$, the eigenvalues of the drifting Hodge Laplacian $\Delta_f$ acting on $p$-forms on a domain $\Omega$ of $M$ with Dirichlet boundary conditions satisfy for any $k\geq 1$
\begin{eqnarray*}
\sum_{i=1}^k(\lambda_{k+1,p,f}-\lambda_{i,p,f})^\alpha&\leq& 
\frac{4}{n}\sum_{i=1}^k(\lambda_{k+1,p,f}-\lambda_{i,p,f})^{\alpha-1}\big(\lambda_{i,p,f}-\delta_1+\frac{1}{4}\delta_2)
\end{eqnarray*}
for $\alpha\leq 2$. Also, we have 
\begin{eqnarray*}
\sum_{i=1}^k(\lambda_{k+1,p,f}-\lambda_{i,p,f})^\alpha&\leq& 
\frac{2\alpha}{n}\sum_{i=1}^k(\lambda_{k+1,p,f}-\lambda_{i,p,f})^{\alpha-1}
\big(\lambda_{i,p,f}-\delta_1+\frac{1}{4}\delta_2)
\end{eqnarray*}
for $\alpha\geq 2$. Here, we set $\delta_2=\mathop{\rm sup}\limits_{\Omega} (n^2|H|^2-2\Delta f-|df|^2)$. 
\end{Cor}

\begin{Bem}
\begin{enumerate}
 \item Choosing \(\alpha=2\) in the estimates of Corollary \ref{cor:recursionformulas} generalizes
the results obtained by Xia and Xu \cite{XX} for functions to the case of \(p\)-forms.
\item In the case of \(\alpha=2\) and \(f=\frac{|X|^2}{2}\) 
recursion formulas similar to the ones of Corollary \ref{cor:recursionformulas} were 
established in \cite{CZ}.
\end{enumerate}
\end{Bem}

Finally, we finish by establishing an estimate for the index of a compact $f$-minimal hypersurface of 
$\R^{n+1}$ as in \cite{S2, IRS}. Recall that a $f$-minimal hypersurface $M$ of the weighted manifold $(\R^{n+1},{\rm can}, e^{-f}dv_g)$ is a hypersurface such that the $f$-mean curvature vanishes,
\begin{align*}
nH_f:=nH+\frac{\partial f}{\partial \nu}=0, 
\end{align*}
where $H$ is the mean curvature. The \emph{stability operator} of $M$, 
also called \emph{Jacobi operator},
is then defined for any smooth function $u$ on $M$ by 
\begin{align*}
L_f(u):=\Delta_f u-({\rm Ric}_f^{\R^{n+1}}(\nu,\nu)+|{\bf II}|^2)u, 
\end{align*}
where \({\rm Ric}_f^{\R^{n+1}}\) denotes the Bakry-Emery Ricci tensor.
On an arbitrary Riemannian manifold  \((N,h)\), the Bakry-Emery Ricci tensor (or the $\infty$-Bakry-Emery Ricci tensor) is defined as follows
\begin{align*}
\operatorname{Ric}^N_f=\operatorname{Ric}^N+\operatorname{Hess}^N f,
\end{align*}
where ${\rm Hess}^{N} f$ is the Hessian of $f$ on $N$.
In the case of \(N=(\R^n,\rm can)\) we simply have that
\({\rm Ric}_f^{\R^{n+1}}={\rm Hess}^{\R^{n+1}} f\)
and the stability operator acquires the form 
\begin{align*}
L_f(u):=\Delta_f u-({\rm Hess}^{\R^{n+1}} f(\nu,\nu)+|{\bf II}|^2)u.
\end{align*}

The $f$-index of $M$ is the number of negative eigenvalues of the Jacobi operator $L_f$. In order to estimate this number, we will establish an upper bound for the eigenvalues 
\(\big(\lambda_k(L_f)\big)_k\) of the Jacobi operator that depends on the eigenvalues of the drifting Hodge Laplacian $\Delta_f$ on $p$-forms. This will generalize the result obtained in \cite[Theorem A]{IRS} when $p=1$. 

\begin{Thm} 
\label{thm:jacobi}
Let $(M^n,g)$ be a compact hypersurface of the weighted manifold $(\R^{n+1}, {\rm can}, d\mu_f=e^{-f}dv_g)$. Assume that $M$ is $f$-minimal and that ${\rm Ric}_f^{\R^{n+1}}={\rm Hess}^{\R^{n+1}} f\geq a>0$. We denote by $k_1,\ldots,k_n$ the principal curvatures of the hypersurface. Then, for all $l\geq 1$ and $1\leq p\leq n$, we have that 
$$\lambda_l(L_f)\leq \lambda_{d(l),p,f}-(p+1)a+\gamma_M p(p-1),$$
where $\gamma_M={\rm sup}\{k_ik_j: i\neq j\}$ and  
$$d(l)=\begin{pmatrix}n+1\\p+1\end{pmatrix}(l-1)+1.$$
\end{Thm}

Now, if we set 
$$\beta=\sharp\{\text{eigenvalues on}\,\, p\text{-forms of}\, \Delta_f\, \text{which are less than}\, (p+1)a-\gamma_M p(p-1)\},$$ 

we obtain the immediate two corollaries 
\begin{Cor} 
\label{cor:jacobi}
Let $(M^n,g)$ be a compact hypersurface of the weighted manifold $(\R^{n+1}, {\rm can}, d\mu_f=e^{-f}dv_g)$. Assume that $M$ is $f$-minimal and that ${\rm Ric}_f^{\R^{n+1}}={\rm Hess}^{\R^{n+1}} f\geq a>0$. Then, for any $1\leq p\leq n$, we have that
$${\rm Ind}_f(M)\geq \frac{1}{\begin{pmatrix}n+1\\p+1\end{pmatrix}}\beta.$$
\end{Cor}

A particular case of Corollary \ref{cor:jacobi} is when $M$ is self-shrinker. Recall that a self-shrinker manifold is a compact and connected  submanifold of the Euclidean space $\R^{n+m}$ satisfying the equation $X^\perp=-nH$. Here $X=(x_1,x_2,\ldots,x_{n+m})$ are the components of the immersion into $\R^{n+m}$. By taking the weight $f=\frac{|X|^2}{2}$, one can easily check that any self-shrinker hypersurface is automatically $f$-minimal and that ${\rm Hess}^{\R^{n+1}} f={\rm can}$. Therefore,
\begin{Cor} \label{cor:jacobibis}
Let $(M^n,g)\to \R^{n+1}$ be a compact self-shrinker such that $p(p-1)\gamma_M\leq p+1$ for some $1\leq p\leq n$. Then, we have 
$${\rm Ind}_f(M)\geq \frac{1}{\begin{pmatrix}n+1\\p+1\end{pmatrix}}b_p(M)+n+1,$$
for the weight $f=\frac{|X|^2}{2}$.
\end{Cor}

\begin{Bem}
We would like to point out that Theorem \ref{thm:jacobi} allows to obtain a statement
on the stability of \(f\)-minimal hypersurfaces in terms of its \(p\)-th Betti number
which was only known for \(p=1\) so far. However, the drawback is the fact that the stability
estimate also depends on the curvature of the hypersurface which is not the case for \(p=1\).
\end{Bem}

This article is organized as follows:
In Section 2 we present a number of general features on Dirac and Laplace type operators
on weighted manifolds. Section 3 provides the proofs of the main results while Section 4 is devoted
to a number of geometric applications. 
Throughout this article, we will use the geometers sign convention for the Laplacian.

\par\medskip
When finalizing this manuscript, the article ``Spinors and mass on weighted manifolds'' \cite{BO}
appeared, where most of the results presented in Subsection 2.1 have been obtained 
independently from our calculations.

\section{Dirac and Laplace type operators on weighted manifolds}

\subsection{The Dirac operator on weighted manifolds}
In this section, we study the question of defining the fundamental Dirac operator
on a weighted Riemannian spin manifold. In the famous article of Perelman \cite[Remark 1.3]{Per}
some results in this direction were already presented.

In particular, this will allow to get a new vanishing result on the kernel of the Dirac operator based on the scalar curvature of weighted manifolds, see also \cite[Remark 1.3]{Per}.
\par\medskip
In order to approach this question, we first recall some basic facts on spin manifolds \cite{LM}. 
Let \((M^n,g)\)
be a Riemannian spin manifold of dimension $n$. The spinor bundle $\Sigma M$ is a vector bundle over the manifold \(M\) that is equipped with
a metric connection \(\nabla\) and a Hermitian scalar product \(\langle\cdot,\cdot\rangle\).
On this bundle we have the Clifford multiplication
of spinors with tangent vectors which we denote by the symbol $``\cdot"$. Remember that Clifford multiplication is compatible with the connection on the spinor bundle. 
Moreover, it is skew-symmetric in the sense
\begin{align*}
\langle Y\cdot\psi,\varphi\rangle
=-\langle \psi,Y\cdot\varphi\rangle
\end{align*}
for all \(Y\in TM\) and \(\psi,\varphi\in\Gamma(\Sigma M)\).
The fundamental Dirac operator \(D\colon\Gamma(\Sigma M)\to\Gamma(\Sigma M)\) acting on spinors is defined by
\begin{align*}
D:=\sum_{i=1}^ne_i\cdot\nabla_{e_i},
\end{align*}
where $\{e_i\}_{i=1,\ldots,n}$ is an orthonormal basis of \(TM\).
The Dirac operator is an elliptic, first order differential operator which is self-adjoint
with respect to the \(L^2\)-norm (when $M$ is compact). Therefore, it admits a sequence of eigenvalues of finite multiplicities. Also, its square is related to the Laplacian through the so-called \emph{Schr\"odinger-Lichnerowicz formula} which is given by \cite{LM}
\begin{align}
\label{eq:schroedinger}
D^2=\nabla^\ast\nabla+\frac{1}{4}\operatorname{Scal}^M,
\end{align}
where $\nabla^*\nabla=-\sum_{i=1}^n \nabla_{e_i}\nabla_{e_i}+\sum_{i=1}^n \nabla_{\nabla_{e_i} e_i}$ and ${\rm Scal}^M$ is the scalar curvature of the manifold $M$. Let us now consider a smooth real valued function $f$ on $M$ and the weighted measure $d\mu_f:=e^{-f}dv_g$.  
   In order to define the Dirac operator on a compact weighted manifold, we use the following variational principle.
The standard Dirac action with a weighted measure is
\begin{align}
\label{action-dirac}
E(\psi):=\int_M\langle\psi,D\psi\rangle d\mu_f,
\end{align}
where \(\psi\in\Gamma(\Sigma M)\).

\begin{Prop}
The critical points of \eqref{action-dirac} are characterized by the equation
\begin{align}
D\psi-\frac{1}{2}df\cdot\psi=0.
\end{align}
\end{Prop}
\begin{proof}
We consider a variation of the spinor \(\psi\) defined as \(\psi_t\colon (-\epsilon,\epsilon)\times M\to \Sigma M\) satisfying \(\frac{\nabla\psi_t}{\partial t}\big|_{t=0}=\varphi\) and we calculate
\begin{align*}
\frac{d}{dt}\big|_{t=0}E(\psi_t)=&\int_M(\langle\varphi,D\psi\rangle +\langle\psi,D\varphi\rangle) d\mu_f \\
=&\int_M\big(\langle\varphi,D\psi\rangle+\langle D\psi,\varphi\rangle-\langle df\cdot\psi,\varphi\rangle\big) d\mu_f\\
=&\int_M\big(\operatorname{Re}\langle (2D-df\cdot)\psi,\varphi\rangle d\mu_f.
\end{align*}
Therefore, the critical points of the energy $E$ are exactly solutions of Equation \eqref{action-dirac}. This completes the proof.
\end{proof}

Motivated by the above calculation, we define
\begin{Dfn}
The fundamental Dirac operator on a Bakry-Emery manifold $(M^n,g, d\mu_f)$ is given by
\begin{align}
\label{dfn:diracbakry}
D_f:=D-\frac{1}{2}df\cdot
\end{align}
\end{Dfn}

\begin{Prop}
The operator \(D_f\) is elliptic. Moreover, if \(M\) is compact, it is self-adjoint with respect
to the \(L^2\)-norm.
\end{Prop}
\begin{proof}
It is straightforward to check that \(D\) and \(D_f\) have the same principal symbol,
hence \(D_f\) is clearly elliptic. To prove the second part, we choose two spinors \(\psi,\varphi\) and calculate
\begin{align*}
\int_M\langle D_f\psi,\varphi\rangle e^{-f} dv_g=&
\int_M\langle D\psi,\varphi\rangle e^{-f} dv_g
-\frac{1}{2}\int_M\langle df\cdot\psi,\varphi\rangle e^{-f} dv_g\\ 
=&\int_M\langle \psi,D(e^{-f}\varphi)\rangle  dv_g
-\frac{1}{2}\int_M\langle df\cdot\psi,\varphi\rangle e^{-f} dv_g\\
=&\int_M\langle\psi,D\varphi\rangle e^{-f} dv_g
-\int_M\langle\psi,df\cdot\varphi\rangle e^{-f} dv_g
+\frac{1}{2}\int_M\langle\psi, df\cdot\varphi\rangle e^{-f} dv_g\\
=&\int_M\langle\psi,D\varphi-\frac{1}{2}df\cdot\varphi\rangle e^{-f} dv_g\\
=&\int_M\langle\psi,D_f\varphi\rangle e^{-f} dv_g.
\end{align*}
The proof is complete.
\end{proof}

\begin{Bem}
It is straightforward to check that \eqref{action-dirac} is invariant under conformal transformations
of the metric on \(M\). Hence, we cannot expect that we can relate the properties of \(D\) and \(D_f\)
by conformally transforming the metric on \(M\).
\end{Bem}

However, the spectral properties of \(D_f\) are essentially the same as the ones of the
standard Dirac operator \(D\) as is shown by the following proposition.

\begin{Prop}\label{prop:invariance}
On a Bakry-Emery manifold $(M^n,g, d\mu_f)$, we have $$D_f(e^{f/2})=e^{f/2} D.$$
In particular, if $\psi$ is an eigenspinor of $D$ associated with the eigenvalue $\lambda$,  then $e^{f/2} \psi$ is an eigenspinor of $D_f$ associated with the same eigenvalue $\lambda$. 
\end{Prop}
\begin{proof}
By a straightforward calculation, we have for any spinor field $\psi$ 
\begin{eqnarray*}
D_f(e^{f/2}\psi)&=&D(e^{f/2}\psi)-\frac{1}{2}e^{f/2}df\cdot\psi\\
&=&e^{f/2}D\psi+\frac{1}{2}e^{f/2}df\cdot\psi-\frac{1}{2}e^{f/2}df\cdot\psi \\
&=&e^{f/2} D\psi.
\end{eqnarray*}
\end{proof}

In order to obtain a corresponding \emph{Schrödinger-Lichnerowicz formula} for \(D_f\)
we make the following observation:
The \(L^2\)-adjoint of \(\nabla\) with respect to the weighted
measure \(e^{-f}dv_g\) will be denoted by \(\nabla^\ast_f\) and is given
by the following expression
\begin{align*}
\nabla^\ast_f=\nabla^\ast+\nabla_{df}.
\end{align*}

Hence, we have 

\begin{Prop}\label{prop:diracbesa}
The square of $D_f^2$ satisfies the following formula: 
\begin{align}
\label{eq:schroedingerbe}
D_f^2=\nabla_f^\ast\nabla+\frac{1}{4}\operatorname{Scal}^M
-\frac{1}{2}\Delta f-\frac{1}{4} |df|^2.
\end{align}
\end{Prop}
\begin{proof}
A direct calculation shows 
\begin{eqnarray*} 
D_f^2&=&(D-\frac{1}{2} df\cdot)(D-\frac{1}{2} df\cdot)\\
&=& D^2-\frac{1}{2}D(df\cdot)-\frac{1}{2}df\cdot D+\frac{1}{4} df\cdot df\cdot\\
&=&D^2-\frac{1}{2}(\delta df-df\cdot D-2\nabla_{df})-\frac{1}{2}df\cdot D+\frac{1}{4} df\cdot df\cdot\\
&=&D^2-\frac{1}{2}\Delta f+\nabla_{df}-\frac{1}{4} |df|^2.
\end{eqnarray*}
Combining with \eqref{eq:schroedinger} and using the definition of \(\nabla^\ast_f\) completes the proof.
\end{proof}

As in the case of the fundamental Dirac operator on Riemannian spin manifolds \cite{Li}, we deduce the following vanishing result  
\begin{Cor}\label{cor:rigidity}
Let \((M,g,d\mu_f)\) be a closed Bakry-Emery manifold. If $\operatorname{Scal}^M>2\Delta f+|df|^2$, then \(\ker D=0\).
\end{Cor}
\begin{proof}
This follows directly from Proposition \ref{prop:diracbesa} and the invariance of the spectrum of the Dirac operator in Proposition \ref{prop:invariance}.
\end{proof}

\begin{Bem}
Suppose that \((M,g)\) is a closed Riemannian surface. If we integrate the condition 
$\operatorname{Scal}^M>2\Delta f+|df|^2$ with respect to the standard Riemannian measure \(dv_g\), we get
\begin{align*}
2\pi\chi(M)-\int_M|df|^2 dv_g>0,
\end{align*}
where \(\chi(M)\) represents the Euler characteristic of the surface.
Hence, this inequality can only be satisfied on a surface of positive Euler characteristic.  
\end{Bem}
 
\subsection{The Bochner-Laplacian on weighted manifolds}  \label{sec:laplacweighted}
In this section, we will recall the Hodge Laplacian on weighted manifolds and define 
a twisted Laplacian motivated from the expression of the Dirac operator in the previous section. This will also allow to get a new vanishing result on the cohomology groups of the manifold.\\

For this, let \((M^n,g,d\mu_f=e^{-f}dv_g)\) be a Riemannian Bakry-Emery manifold. We denote by
$d$ the exterior differential and by $\delta$ the codifferential. The weighted codifferential is defined by $\delta_f:=\delta+df\lrcorner$  where we identify here (and in all the paper) vectors with one-forms through the musical isomorphism. The weighted codifferential is the $L^2$-formal adjoint of $d$ with respect to the measure $d\mu_f$, when $M$ is compact. By a straightforward computation, one shows that 
$$\delta_f^2=\delta(\delta+df\lrcorner)+df\lrcorner(\delta+df\lrcorner)=\delta(df\lrcorner)+df\lrcorner\delta=0,$$
where we use the fact $\delta(df\lrcorner)=-df\lrcorner\delta$, since ${\rm Hess}^Mf=\nabla df$ is a symmetric endomorphism.  The drifting Hodge Laplacian on differential forms is then defined as $$\Delta_f:=d\delta_f+\delta_f d.$$ 
Clearly, the drifting Hodge Laplacian commutes with $d$ and $\delta_f$ and, therefore, it preserves the spaces of exact and weighted coexact forms. Moreover, this operator is elliptic and, when $M$ is compact, it is self-adjoint with respect to the weighted measure $d\mu_f$. Therefore, as for the ordinary Hodge Laplacian on compact manifolds, the drifting Hodge Laplacian restricted to differential $p$-forms ($1\leq p\leq n$) has a spectrum that consists of a nondecreasing, unbounded sequence of eigenvalues  with finite multiplicities, that is
\begin{align*}
\operatorname{Spec}(\Delta_f)=\{0\leq\lambda_{1,p,f}\leq\lambda_{2,p,f}\leq\ldots\}. 
\end{align*} 
The eigenvalue $0$ corresponds to the space of $f$-harmonic forms, that is differential forms $\omega$  satisfying $d\omega=0$ and $\delta_f\omega=0$. Notice that, by standard elliptic theory,  we have an isomorphism \cite[Formula 2.13]{L},  
$$H^p(M)\simeq\{\omega\in \Omega^p(M)|\,\, d\omega=0,\, \delta_f\omega=0\},$$ 
meaning that the $f$-Betti numbers do not depend on $f$.  Now it is shown in \cite{PW}, that the drifting Hodge Laplacian $\Delta_f$ has a corresponding Bochner-Weitzenb\"ock formula which is 
\begin{equation}\label{eq:Bochner}
\Delta_f=\nabla_f^*\nabla+\mathfrak{B}^{[p]}+T_f^{[p]}
\end{equation}
where $\mathfrak{B}^{[p]}=\sum_{i,j=1}^ne_j^*\wedge e_i\lrcorner R^M(e_i,e_j)$ is the Bochner operator that appears in the Bochner-Weitzenb\"ock formula for $\Delta=d\delta+\delta d$. Here, $R^M(X,Y):=[\nabla_X,\nabla_Y]-\nabla_{[X,Y]}$ is the curvature tensor operator of $M$ for $X,Y\in TM$ and $T_f^{[p]}$ is the self-adjoint endomorphism of $\Lambda^p(M)$ given by 
\begin{align}
(T_f^{[p]}\omega)(X_1,\ldots,X_p)=\sum_{k=1}^p\omega(X_1,\ldots,\nabla^2f(X_k),\ldots,X_p)
\end{align}
for all $X_1,\ldots, X_p\in TM$ and $\{e_1,\ldots,e_n\}$ is a local orthonormal frame of $TM$. The tensor $\mathfrak{B}^{[p]}+T_f^{[p]}$ is called the \emph{$p$-Ricci tensor}
and is denoted by ${\rm Ric}^{(p)}_f$ (see \cite{P}). For $p=1$, the $p$-Ricci tensor is just ${\rm Ric}^{(1)}_f={\rm Ric}^M+{\rm Hess}^M f$ which is the $\infty$-Bakry-Emery Ricci tensor.
Also, for any $N\in ]-\infty,0[\cup ]n-p+1,\infty[$, 
we let $${\rm Ric}^{(p)}_{N,f}:= {\rm Ric}^{(p)}_{f}-\frac{1}{N-(n-p+1)}(df\wedge (df\lrcorner))$$ which corresponds to the so-called $N$-Bakry-Emery Ricci tensor for $p=1$.  It is now a well-known fact that the Bochner-Weitzenb\"ock formula gives rise to the Gallot-Meyer estimate on manifolds having a lower bound on the Bochner operator $\mathfrak{B}^{[p]}$ \cite{GM}. In the following, we will adapt this technique to give an estimate for the eigenvalues of the drifting Hodge Laplacian on a Bakry-Emery manifold having a lower bound on ${\rm Ric}^{(p)}_{N,f}$. For this, we denote by $\lambda'_{1,p,f}$ the first positive eigenvalue of $\Delta_f$ restricted to exact $p$-forms. Then, we have 

\begin{Prop} Let $(M^n,g,d\mu_f)$ be a compact Bakry-Emery manifold. If ${\rm Ric}^{(p)}_{N,f}\geq p(n-p)\gamma$ for some $\gamma>0$, then the first eigenvalue $\lambda'_{1,p,f}$ satisfies the estimate 
$$\lambda'_{1,p,f}\geq p(n-p)\gamma\frac{N}{N-1}.$$
\end{Prop}

\begin{proof} 
Let $\omega$ be any exact $p$-eigenform. Applying the Bochner-Weitzenb\"ock formula \eqref{eq:Bochner} yields 
$$\lambda'_{1,p,f}\int_M|\omega|^2d\mu_f=\int_M|\nabla\omega|^2d\mu_f+\int_M\langle(\mathfrak{B}^{[p]}+T_f^{[p]})\omega,\omega\rangle d\mu_f.$$ 
Now, using that $|\nabla\omega|^2\geq \frac{1}{n-p+1}|\delta\omega|^2$ as $\omega$ is closed \cite{GM}, we have that 
\begin{eqnarray*}
|\nabla\omega|^2&\geq&\frac{1}{n-p+1}|\delta_f\omega-df\lrcorner\omega|^2\\
&\geq &\frac{1}{n-p+1}(|\delta_f\omega|-|df\lrcorner\omega|)^2\\
&\geq &\frac{1}{N}|\delta_f\omega|^2-\frac{1}{N-(n-p+1)}|df\lrcorner\omega|^2. 
\end{eqnarray*}
Here, we use the inequality $\frac{(p+q)^2}{s}\geq \frac{p^2}{N}-\frac{q^2}{N-s}$ for all $N$ such that $N(N-s)>0$. Therefore, by integrating over $M$ and using the fact $\lambda'_{1,p,f}\int_M|\omega|^2d\mu_f=\int_M|\delta_f\omega|^2d\mu_f$, we get the estimate.
\end{proof}

In the last part of this section, we will define a new drifting Hodge Laplacian acting on differential $p$-forms that will allow us to get a new vanishing result on the cohomology groups. For this, recall that the weighted Dirac operator defined in the previous section is $D_f=\sum_{i=1}^n e_i\cdot\nabla_{e_i}-\frac{1}{2}df\cdot$. When restricted to differential forms, the operator $D_f$ can be written as 
$$D_f=\sum_{i=1}^n e_i\wedge\nabla_{e_i}-e_i\lrcorner\nabla_{e_i}-\frac{1}{2}df\wedge+\frac{1}{2}df\lrcorner=\widetilde{d}_f+\widetilde{\delta}_f,$$ 
where $\widetilde{d}_f:=d-\frac{1}{2}df\wedge$ and  
$\widetilde{\delta}_f:=\delta+\frac{1}{2}df\lrcorner$. Here, we use the fact that $X\cdot\omega=X\wedge\omega-X\lrcorner\omega$ for any vector field $X$ and a differential form $\omega$.  It is not difficult to check that $\widetilde{d}_f^2=\widetilde{\delta}_f^2=0$ and that $\widetilde{\delta}_f=\delta_f-\frac{1}{2}df\lrcorner$ is the $L^2$-adjoint of $\widetilde{d}_f$ with respect to the measure $d\mu_f=e^{-f}dv_g$.  Also the square of the Dirac operator gives rise to the twisted Laplacian 
$$D_f^2=\widetilde{ \Delta}_f:=\widetilde{d}_f\widetilde{\delta}_f+\widetilde{\delta}_f\widetilde{d}_f$$ 
that has the same spectrum as $\Delta$ by Proposition \ref{prop:invariance}. Now, we establish a Bochner-Weitzenb\"ock formula for  $\widetilde{\Delta}_f$ in order to get a new vanishing result for the cohomology groups.

\begin{Prop} Let $(M,g,d\mu_f)$ be a Bakry-Emery manifold. Then, we have the Bochner-Weitzenb\"ock formula for the twisted drifting Hodge Laplacian 
\begin{equation}\label{eq:bochnertildelaplace}
\widetilde{ \Delta}_f=\nabla_f^*\nabla+\mathfrak{B}^{[p]}-\frac{1}{2}\Delta f-\frac{1}{4}|df|^2.
\end{equation}
In particular, if $M$ is compact and $\mathfrak{B}^{[p]}>\frac{1}{2}\Delta f+\frac{1}{4}|df|^2$ for some $p$, then $H^p(M)=0$.
\end{Prop}
\begin{proof} We compute
\begin{eqnarray*} 
\widetilde{ \Delta}_f&=&\widetilde{d}_f\widetilde{\delta}_f+\widetilde{\delta}_f\widetilde{d}_f\\
&=& (d-\frac{1}{2}df\wedge)(\delta_f-\frac{1}{2} df\lrcorner)+(\delta_f-\frac{1}{2} df\lrcorner)(d-\frac{1}{2}df\wedge)\\
&=&\Delta_f-\frac{1}{2}d(df\lrcorner)-\frac{1}{2}df\wedge\delta_f+\frac{1}{4}df\wedge (df\lrcorner)-\frac{1}{2}\delta_f(df\wedge)-\frac{1}{2}df\lrcorner d+\frac{1}{4}df\lrcorner (df\wedge)\\
&=&\Delta_f-\frac{1}{2}\mathcal{L}_{df}-\frac{1}{2}df\wedge\delta_f-\frac{1}{2}\delta_f(df\wedge)+\frac{1}{4}|df|^2\\
&=&\Delta_f-\frac{1}{2}\nabla_{df}-\frac{1}{2}T^{[p]}_f-\frac{1}{2}df\wedge\delta_f-\frac{1}{2}\delta_f(df\wedge)+\frac{1}{4}|df|^2.
\end{eqnarray*}
In the last equality, we used the identity \cite[Lem. 2.1]{S1}
$$\mathcal{L}_{df}=\nabla_{df}+T_f^{[p]}.$$ 
Now, an easy computation shows that 
$$\delta_f(df\wedge)=\Delta f+\sum_{i=1}^n\nabla_{e_i}df\wedge (e_i\lrcorner)-\nabla_{df}+|df|^2-df\wedge\delta_f.$$
Using the fact that $T_f^{[p]}=\sum_{i=1}^n\nabla_{e_i}df\wedge (e_i\lrcorner)$ which can be proven by a straightforward computation,  we get after using Equation \eqref{eq:Bochner} and replacing the last equality,
\begin{eqnarray*}
\widetilde{ \Delta}_f&=& \nabla_f^*\nabla+\mathfrak{B}^{[p]}+T_f^{[p]}-\frac{1}{2}\nabla_{df}-\frac{1}{2}T^{[p]}_f-\frac{1}{2}df\wedge\delta_f \\
\nonumber&&-\frac{1}{2}\left(\Delta f+T_f^{[p]}-\nabla_{df}+|df|^2-df\wedge\delta_f\right)
+\frac{1}{4}|df|^2\\
\nonumber&=&\nabla_f^*\nabla+\mathfrak{B}^{[p]}-\frac{1}{2}\Delta f-\frac{1}{4}|df|^2.
\end{eqnarray*}
 Note that \eqref{eq:bochnertildelaplace} has the same structure as the corresponding
Schr\"odinger-Lichnerowicz formula for the Dirac operator on weighted manifolds \eqref{eq:schroedingerbe}. The vanishing result on the cohomology is obtained by just applying  \eqref{eq:bochnertildelaplace} to a $\widetilde{ \Delta}_f$-harmonic form and integrating over $M$. Finally, the fact that the set of $\widetilde{ \Delta}_f$-harmonic form is isomorphic to $H^p(M)$ allows to finish the proof.
\end{proof}

\section{Proofs of the main results}
In this section, we provide the proofs of the main results.
\begin{proof}[Proof of Theorem \ref{thm:inequality}]  
We start by proving the first part of the theorem. We assume that $(M^n,g, d\mu_f=e^{-f}dv_g)$ is a compact Bakry-Emery manifold that is isometrically immersed into $\R^{n+m}$. For $p=1,\ldots,n$, we let $\lambda_{1,p,f}$ the first nonegative eigenvalue of $\Delta_f$ on $p$-forms and $\lambda'_{1,p,f}$ the first positive eigenvalue of $\Delta_f$ on exact $p$-forms.  For every $i=1,\ldots,n+m$, the parallel unit vector field $\partial x_i$ decomposes into $\partial x_i=(\partial x_i)^T+(\partial x_i)^\perp$ with $(\partial x_i)^T=d(x_i\circ \iota)$ where $\iota$ is the isometric immersion.  
Now, take any $p$-eigenform $\omega$ of $\Delta_f$ and consider the $(p-1)$-form $\phi_i:=(\partial x_i)^T\lrcorner \omega$ for each $i=1,\ldots, n+m$. 
By the Rayleigh min-max principle, we have that 
\begin{equation}\label{eq:minmaxdrifting}
\lambda_{1,p-1,f}\sum_{i=1}^{n+m}\int_M|\phi_i|^2d\mu_f\leq \sum_{i=1}^{n+m}\int_M(|d\phi_i|^2
+|\delta_f\phi_i|^2)d\mu_f.
\end{equation}
In the following, we will adapt some computations done in \cite{GS} to the context of the drifting Hodge Laplacian (see also \cite{Asa}, \cite[Thm. 5.8]{CGH} for a similar computation). 
Denoting by $\{e_1,\ldots, e_n\}$ a local orthonormal frame of $TM$
we find that
\begin{align}
\label{identity-a}
\sum_{i=1}^{n+m}|\phi_i|^2=\sum_{s,t=1}^n\underbrace{\sum_{i=1}^{n+m}g((\partial x_i)^T,e_s)g((\partial x_i)^T,e_t)}_{\delta_{st}}\langle e_s\lrcorner \omega,e_t\lrcorner\omega\rangle=\sum_{s=1}^n |e_s\lrcorner \omega|^2=p|\omega|^2. 
\end{align}
Here, we use the fact that $\alpha=\frac{1}{p}\sum_{s=1}^n e_s^*\wedge e_s\lrcorner\alpha$, for any $p$-form $\alpha$.
Moreover, by using the Cartan formula and \cite[Eq. (4.3)]{GS}, we write
$$d\phi_i=\mathcal{L}_{(\partial x_i)^T}\omega-(\partial x_i)^T\lrcorner d\omega=\nabla_{(\partial x_i)^T}\omega+{\bf II}_{(\partial x_i)^\perp}^{[p]}\omega-(\partial x_i)^T\lrcorner d\omega,$$
where ${\bf {\bf II}}_{Z}^{[p]}$ is the canonical extension of the second fundamental form ${\bf II}$ of the immersion to differential $p$-forms in the normal direction $Z$. More precisely, if we write $\langle {\bf {\bf II}}_Z(X),Y\rangle=\langle {\bf {\bf II}}(X,Y),Z\rangle$ for any tangent vector fields $X,Y\in TM$ and $Z\in T^\perp M$, we define 
\begin{align}
\label{dfn:sff-forms}
({\bf {\bf II}}_Z^{[p]}\omega)(X_1,\ldots, X_p):=\sum_{i=1}^p \omega(X_1,\ldots, {\bf {\bf II}}_Z(X_i),\ldots, X_p)
\end{align}
for any differential $p$-form $\omega$ on $M$ and \(X_1,\ldots,X_p\in TM\). Now, as we did in \eqref{eq:minmaxdrifting}, we decompose $(\partial x_i)^T$ and $(\partial x_i)^\perp$ in the frames $\{e_1,\ldots,e_n\}$ and $\{\nu_1,\ldots,\nu_m\}$ to compute  the norm square of $d\phi_i$, we deduce after summing over $i$ that 
\begin{eqnarray}\label{eq:dfomega}
\sum_{i=1}^{n+m}|d\phi_i|^2&=&|\nabla\omega|^2+\sum_{s=1}^m|{\bf {\bf II}}^{[p]}_{\nu_s}\omega|^2+(p+1)|d\omega|^2-2\sum_{i=1}^n \langle\nabla_{e_i}\omega,e_i\lrcorner d\omega\rangle\nonumber\\
&=&|\nabla\omega|^2+\sum_{s=1}^m|{\bf {\bf II}}^{[p]}_{\nu_s}\omega|^2+(p-1)|d\omega|^2.
\end{eqnarray}
Here, $\{\nu_1,\ldots,\nu_m\}$ is a local orthonormal frame of $T^\perp M$. In the above computation, we use the fact that all cross terms involving $(\partial x_i)^T$ and $(\partial x_i)^\perp$ are zero.  
By writing $\delta_f\phi_i=\delta \phi_i+df\lrcorner \phi_i$ with $\delta\phi_i=-(\partial x_i)^T\lrcorner\delta\omega$ which is a consequence from the fact that $\nabla(\partial x_i)^T={\rm Hess}^M(x_i\circ\iota)$ is a symmetric endomorphism, we get that 
\begin{eqnarray}\label{eq:deltafomega}
 \sum_{i=1}^{n+m}|\delta_f\phi_i|^2&=&
(p-1)|\delta\omega|^2+(p-1)|df\lrcorner\omega|^2
-2\sum_{i=1}^{n+m}\langle (\partial x_i)^T\lrcorner\delta\omega, df\lrcorner(\partial x_i)^T\lrcorner\omega\rangle\nonumber\\
 &=&(p-1)|\delta\omega|^2+(p-1)|df\lrcorner\omega|^2
 -2\sum_{i=1}^{n}\langle e_i\lrcorner \delta\omega, df\lrcorner e_i\lrcorner\omega\rangle\nonumber\\
 &=&(p-1)|\delta\omega|^2+(p-1)|df\lrcorner\omega|^2+2(p-1)\langle\delta\omega,df\lrcorner\omega\rangle\nonumber\\
 &=&(p-1)|\delta_f\omega|^2.
 \end{eqnarray}
Replacing Equalities \eqref{identity-a}, \eqref{eq:dfomega} and \eqref{eq:deltafomega} into Inequality \eqref{eq:minmaxdrifting}, we deduce after using the Bochner-Weitzenb\"ock formula \eqref{eq:Bochner} that 
\begin{eqnarray*}
p\lambda_{1,p-1,f}\int_M|\omega|^2d\mu_f&\leq& \int_M(|\nabla\omega|^2+\sum_{s=1}^m|{\bf {\bf II}}^{[p]}_{\nu_s}\omega|^2+(p-1)|d\omega|^2+(p-1)|\delta_f\omega|^2)d\mu_f\\
&=&\int_M(\langle\Delta_f\omega,\omega\rangle-\langle(\mathfrak{B}^{[p]}+T_f^{[p]})\omega,\omega\rangle 
+\sum_{s=1}^m|{\bf {\bf II}}^{[p]}_{\nu_s}\omega|^2)d\mu_f \\
&&+(p-1)\lambda_{1,p,f}\int_M|\omega|^2d\mu_f\\
&=&p\lambda_{1,p,f}\int_M|\omega|^2d\mu_f+\int_M\langle\left(\sum_{s=1}^m({\bf {\bf II}}^{[p]}_{\nu_s})^2-\mathfrak{B}^{[p]}-T_f^{[p]}\right)\omega,\omega\rangle d\mu_f.
\end{eqnarray*}
This finishes the proof of the first part. 
To prove the second part of Theorem \ref{thm:inequality}, we use the same technique as in \cite{S} and adapt it to the case of the drifting Hodge Laplacian.  
We consider the exact $p$-form 
$$(\partial x_{i_1})^T\wedge \ldots \wedge (\partial x_{i_p})^T,$$
for $i_k=1,\ldots, n+m, k=1,\ldots,p$. By applying the min-max principle and summing over $i_1,\ldots,i_p$, we get 
\begin{equation}\label{eq:estimate}
\lambda'_{1,p,f}\sum_{i_1,\ldots,i_p}\int_M |(\partial x_{i_1})^T\wedge \ldots \wedge (\partial x_{i_p})^T|^2  d\mu_f\leq \sum_{i_1,\ldots,i_p}\int_M|\delta_f\left((\partial x_{i_1})^T\wedge \ldots \wedge (\partial x_{i_p})^T\right)|^2 d\mu_f.
\end{equation}
Now, we manipulate both sums the same way as in \cite{S}. 
The sum on the left hand side of
\eqref{eq:estimate}
is equal to $p!\begin{pmatrix}n\\p\end{pmatrix}$ (see \cite[Lemma 2.1]{S} for more details). To confirm this fact, we 
denote as usual by $\{e_1,\ldots,e_n\}$ a local orthonormal frame of $TM$ and calculate
\begin{align*}
\sum_{i_1,\ldots,i_p}|(\partial x_{i_1})^T&\wedge \ldots \wedge (\partial x_{i_p})^T|^2 \\
&=\sum_{i_1,\ldots,i_p}\langle (\partial x_{i_1})^T\wedge \ldots \wedge (\partial x_{i_p})^T,(\partial x_{i_1})^T\wedge \ldots \wedge (\partial x_{i_p})^T\rangle\\
&=\sum_{i_1,\ldots,i_p,s,t}g((\partial x_{i_1})^T,e_s)g((\partial x_{i_1})^T,e_t)\langle e_s\wedge \ldots \wedge (\partial x_{i_p})^T,e_t\wedge \ldots \wedge (\partial x_{i_p})^T\rangle\\
&=\sum_{i_2,\ldots,i_p,s,t}\delta_{st}\langle e_s\wedge (\partial x_{i_2})^T\wedge \ldots \wedge (\partial x_{i_p})^T,e_t\wedge (\partial x_{i_2})^T\wedge\ldots \wedge (\partial x_{i_p})^T\rangle\\
&=\sum_{i_2,\ldots,i_p,s}|e_s\wedge (\partial x_{i_2})^T\wedge \ldots \wedge (\partial x_{i_p})^T|^2\\
&= (n-p+1)\sum_{i_2,\ldots,i_p}|(\partial x_{i_2})^T\wedge \ldots \wedge (\partial x_{i_p})^T|^2. 
\end{align*}
In the last line, we use the fact that, for any $p$-form $\omega$, the equality 
$\sum_{s=1}^n |e_s\wedge\omega|^2=(n-p)|\omega|^2$ holds true. Hence, by induction, we deduce that 
\begin{equation}\label{eq:wedgenorm}
\sum_{i_1,\ldots,i_p}|(\partial x_{i_1})^T\wedge \ldots 
\wedge (\partial x_{i_p})^T|^2=(n-p+1)(n-p+2)\ldots n=p!
\begin{pmatrix}n\\p\end{pmatrix}.
\end{equation}

Now, we aim to compute the sum on the right hand side of \eqref{eq:estimate}.
To this end we recall some computations done in \cite{S}. First, by \cite[Eq. 3.3]{S}, 
$$\delta\left((\partial x_{i_1})^T\wedge \ldots \wedge (\partial x_{i_p})^T\right)=\sum_{s=1}^p (-1)^{s+1} T^{[p-1]}_{(\partial x_{i_s})^\perp}((\partial x_{i_1})^T\wedge \ldots \wedge \widehat{(\partial x_{i_s})^T}\wedge\ldots\wedge (\partial x_{i_p})^T)$$
where, on differential $q$-forms, the operator $T_{\nu}^{[q]}$ is defined as 
\begin{align*}
T_{\nu}^{[q]}:={\bf {\bf II}}_{\nu}^{[q]}-n \langle H,\nu\rangle I
\end{align*}
for any normal vector field $\nu$. Here, ${\bf {\bf II}}_{\nu}^{[q]}$ is the canonical extension of the second fundamental form as defined previously in \eqref{dfn:sff-forms}.
For any basis of orthonormal vector fields
$\{\nu_1,\ldots,\nu_m\}\in T_xM^\perp$, we set \\
$$||{\bf {\bf II}}^{[q]}||^2=\sum_{s=1}^m ||{\bf {\bf II}}_{\nu_s}^{[q]}||^2 \quad\text{and}\quad ||T^{[q]}||^2=\sum_{s=1}^m ||T_{\nu_s}^{[q]}||^2.$$

It was shown in  \cite[p. 592]{S} that  
\begin{equation}\label{eq:normdelta}
\sum_{i_1,\ldots,i_p}|\delta\left((\partial x_{i_1})^T\wedge \ldots \wedge (\partial x_{i_p})^T\right)|^2=p!||T^{[p-1]}||^2
\end{equation}
and the latter can be computed explicitly in terms of $|{\bf {\bf II}}|^2$ and the scalar curvature of $M$ through the formula (see \cite[Lemma 2.5]{S})
 $$||T^{[p-1]}||^2=\begin{pmatrix}n\\p\end{pmatrix}\left(pn|H|^2-\frac{p(p-1)}{n(n-1)}{\rm Scal}^M\right).$$

Therefore, we compute
\begin{align}\label{eq:deltaf}
|\delta_f\left((\partial x_{i_1})^T\wedge \ldots \wedge (\partial x_{i_p})^T\right)|^2=
&|\delta\left((\partial x_{i_1})^T\wedge \ldots \wedge (\partial x_{i_p})^T\right)|^2+|df\lrcorner \left((\partial x_{i_1})^T\wedge \ldots \wedge (\partial x_{i_p})^T\right)|^2\nonumber\\
&+2\langle \delta\left((\partial x_{i_1})^T\wedge \ldots \wedge (\partial x_{i_p})^T\right), df\lrcorner \left((\partial x_{i_1})^T\wedge \ldots \wedge (\partial x_{i_p})^T\right)\rangle.
\end{align}
First, we show that the mixed term in \eqref{eq:deltaf}  vanishes. For this, we write
\begin{align*}
\langle & \delta\left((\partial x_{i_1})^T\wedge \ldots \wedge (\partial x_{i_p})^T\right), df\lrcorner \left((\partial x_{i_1})^T\wedge \ldots \wedge (\partial x_{i_p})^T\right)\rangle\\
&=\sum_{s=1}^p (-1)^{s+1}\langle  T^{[p-1]}_{(\partial x_{i_s})^\perp}((\partial x_{i_1})^T\wedge \ldots \wedge \widehat{(\partial x_{i_s})^T}\wedge\ldots\wedge (\partial x_{i_p})^T), df\lrcorner \left((\partial x_{i_1})^T\wedge \ldots \wedge (\partial x_{i_p})^T\right)\rangle \\
&=
\sum_{s=1}^p (-1)^{s+1}g((\partial x_{i_s})^\perp,\nu_\alpha) g((\partial x_{i_s})^T,e_\beta) \\
&\hspace{0.5cm}\times\langle  T^{[p-1]}_{\nu_\alpha}((\partial x_{i_1})^T\wedge \ldots \wedge \widehat{(\partial x_{i_s})^T}\wedge\ldots\wedge (\partial x_{i_p})^T), df\lrcorner \left((\partial x_{i_1})^T\wedge \ldots \wedge e_\beta\wedge \ldots \wedge (\partial x_{i_p})^T\right)\rangle \\
&=
\sum_{s=1}^p (-1)^{s+1}g(\partial x_{i_s},\nu_\alpha) g(\partial x_{i_s},e_\beta) \\
&\hspace{0.5cm}\times\langle  T^{[p-1]}_{\nu_\alpha}((\partial x_{i_1})^T\wedge \ldots \wedge \widehat{(\partial x_{i_s})^T}\wedge\ldots\wedge (\partial x_{i_p})^T), df\lrcorner \left((\partial x_{i_1})^T\wedge \ldots \wedge  e_\beta\wedge \ldots \wedge (\partial x_{i_p})^T\right)\rangle.
\end{align*}

Hence, summing over $i_1,\ldots, i_p$, we deduce that 
$$\sum_{i_1,\ldots,i_p}\langle \delta\left((\partial x_{i_1})^T\wedge \ldots \wedge (\partial x_{i_p})^T\right), df\lrcorner \left((\partial x_{i_1})^T\wedge \ldots \wedge (\partial x_{i_p})^T\right)\rangle=0.$$ Now, it remains to compute the sum over $i_1,\ldots, i_p$ of the second term in  \eqref{eq:deltaf}. We establish the following lemma: 

\begin{Lem} \label{Lem:interiorproduct}For any vector field $X\in TM$, we have 
$$\sum_{i_1,\ldots,i_p} |X\lrcorner \left((\partial x_{i_1})^T\wedge \ldots \wedge (\partial x_{i_p})^T\right)|^2=p!\begin{pmatrix}n-1\\p-1\end{pmatrix}|X|^2.$$ 
\end{Lem}

\begin{proof} Using the formula $X\lrcorner (\alpha\wedge\omega)=g(X,\alpha^\sharp)\omega-\alpha\wedge X\lrcorner \omega$, which is valid for any one-form $\alpha$, we write  
\begin{align*}X\lrcorner \left((\partial x_{i_1})^T\wedge \ldots \wedge (\partial x_{i_p})^T\right)
=&g(X,(\partial x_{i_1})^T)(\partial x_{i_2})^T\wedge \ldots \wedge (\partial x_{i_p})^T \\
&-(\partial x_{i_1})^T\wedge \left(X\lrcorner((\partial x_{i_2})^T\ldots \wedge (\partial x_{i_p})^T)\right). 
\end{align*}

Taking the norm and summing over $i_1,\ldots, i_p$, we get that 
\begin{align*}
&\sum_{i_1,\ldots,i_p} |X\lrcorner \left((\partial x_{i_1})^T\wedge \ldots \wedge (\partial x_{i_p})^T\right)|^2\\
=&\sum_{i_1,\ldots,i_p}  g(X,(\partial x_{i_1})^T)^2|(\partial x_{i_2})^T\wedge \ldots \wedge (\partial x_{i_p})^T|^2+\sum_{i_1,\ldots,i_p}|(\partial x_{i_1})^T\wedge \left(X\lrcorner((\partial x_{i_2})^T\wedge\ldots \wedge (\partial x_{i_p})^T)\right)|^2\\
&-2\sum_{i_1,\ldots,i_p}g(X,(\partial x_{i_1})^T)\langle(\partial x_{i_1})^T\lrcorner\left((\partial x_{i_2})^T\wedge \ldots \wedge (\partial x_{i_p})^T\right), X\lrcorner((\partial x_{i_2})^T\wedge\ldots \wedge (\partial x_{i_p})^T)\rangle.
\end{align*}
Using the fact that $\sum_{i=1}^{n+m} g(X,(\partial x_i)^T)^2=|X|^2$ and Equation \eqref{eq:wedgenorm}, the first sum in the above equation is equal to $(p-1)!\begin{pmatrix}n\\p-1\end{pmatrix}|X|^2$. Concerning the second sum,  we make use of the identity $\sum_{i=1}^{n+m} |(\partial x_i)^T\wedge\omega|^2=(n-p)|\omega|^2$, 
and deduce that the second sum is equal to 
$$(n-p+2)\sum_{i_2,\ldots,i_p} |X\lrcorner \left((\partial x_{i_2})^T\wedge \ldots \wedge (\partial x_{i_p})^T\right)|^2.$$
After using  $X=\sum_{i=1}^{n+m} g(X,(\partial x_i)^T) (\partial x_i)^T$, the last sum reduces to 
$$-2\sum_{i_2,\ldots,i_p} |X\lrcorner \left((\partial x_{i_2})^T\wedge \ldots \wedge (\partial x_{i_p})^T\right)|^2.$$
Here, we conclude that
\begin{align*}
\sum_{i_1,\ldots,i_p} |X\lrcorner \left((\partial x_{i_1})^T\wedge \ldots \wedge (\partial x_{i_p})^T\right)|^2=&(p-1)!\begin{pmatrix}n\\p-1\end{pmatrix}|X|^2 \\
&+(n-p)\sum_{i_2,\ldots,i_p} |X\lrcorner \left((\partial x_{i_2})^T\wedge \ldots \wedge (\partial x_{i_p})^T\right)|^2.
\end{align*}
Thus, by induction, we deduce that 
\begin{align}\label{eq:combinatory}
\sum_{i_1,\ldots,i_p}|X\lrcorner \left((\partial x_{i_1})^T\wedge \ldots \wedge (\partial x_{i_p})^T\right)|^2
\nonumber=&|X|^2\big((p-1)!\begin{pmatrix}n\\p-1\end{pmatrix}
+(n-p)(p-2)!\begin{pmatrix}n\\p-2\end{pmatrix} \\
&+(n-p)(n-p+1)(p-3)!\begin{pmatrix}n\\p-3\end{pmatrix}+\ldots\big).
\end{align}

After some algebraic manipulations, Equality \eqref{eq:combinatory} reduces to the following 
\begin{align*}
|X|^2 \frac{n!}{(n-p-1)!}&\sum_{k=0}^{p-1}\frac{1}{(n-p+k)(n-p+k+1)} \\
=&|X|^2 \frac{n!}{(n-p-1)!}\sum_{k=0}^{p-1}\left(\frac{1}{n-p+k}-\frac{1}{n-p+k+1}\right)\\
=&|X|^2 \frac{n!}{(n-p-1)!} \frac{p}{n(n-p)}\\
=&p!\begin{pmatrix}n-1\\p-1\end{pmatrix}|X|^2.
\end{align*}
This finishes the proof of the Lemma.
\end{proof} 

Using Lemma \ref{Lem:interiorproduct} with $X=df$ and Equation \eqref{eq:normdelta}, we deduce that the sum over $i_1,\ldots,i_p$ of Equation \eqref{eq:deltaf} gives the following 
$$\sum_{i_1,\ldots,i_p}|\delta_f\left((\partial x_{i_1})^T\wedge \ldots \wedge (\partial x_{i_p})^T\right)|^2=p!||T^{[p-1]}||^2+p!\begin{pmatrix}n-1\\p-1\end{pmatrix}|df|^2.$$ 
Hence, by plugging Equation \eqref{eq:wedgenorm} into Inequality \eqref{eq:estimate} we 
get the estimate 
$$\lambda'_{1,p,f}\leq \frac{1}{{\rm Vol}_f(M)}\int_M \left(pn|H|^2-\frac{p(p-1)}{n(n-1)}{\rm Scal}^M+\frac{p}{n}|df|^2\right) d\mu_f$$ 
completing the proof of Theorem \ref{thm:inequality}.
\end{proof}

Before we turn to the proof of Theorem \ref{thm:recursionformulas}, 
we recall the following powerful results from \cite{AH} in which 
the authors show the following: 
\begin{Thm} \cite[Thm. 2]{AH} 
Let $\mathcal{H}$ be a complex Hilbert space with a given inner product $\langle\cdot,\cdot\rangle$. Let $A:\mathcal{D}\subset \mathcal{H}\to \mathcal{H}$ be a self-adjoint operator defined on a dense domain $\mathcal{D}$ which is semi-bounded below and has a discrete spectrum $\lambda_1\leq \lambda_2\leq \lambda_3\ldots$. Let $\{B_k:A(\mathcal{D})\to \mathcal{H}\}_{k=1}^N$ be a collection of symmetric operators which leave $\mathcal{D}$ invariant and let $\{u_i\}_{i=1}^\infty$ be the normalized eigenvectors of $A$, $u_i$ corresponding to $\lambda_i$. This family is assumed to be an orthonormal basis of $\mathcal{H}$. 
Let $g$ be a nonnegative and nondecreasing function of the eigenvalues $\{\lambda_i\}_{i=1}^m$. Then we have the inequality 
$$\sum_{i=1}^m\sum_{k=1}^N (\lambda_{m+1}-\lambda_i)^2g(\lambda_i)\langle [A,B_k]u_i,B_k u_i\rangle\leq \sum_{i=1}^m\sum_{k=1}^N (\lambda_{m+1}-\lambda_i)g(\lambda_i)||[A,B_k]u_i||^2.$$
Here $[A,B]:=AB-BA$ is the commutator of the two operators $A$ and $B$. 
\end{Thm}

As a corollary of this result and by taking $g(\lambda)=(\lambda_{m+1}-\lambda)^{\alpha-2}$ for $\alpha\leq 2$, we get the following inequality \cite[Cor. 3]{AH}
\begin{equation}\label{eq:AH1}
\sum_{i=1}^m\sum_{k=1}^N(\lambda_{m+1}-\lambda_i)^\alpha \langle [A,B_k]u_i,B_k u_i\rangle\leq \sum_{i=1}^m\sum_{k=1}^N (\lambda_{m+1}-\lambda_i)^{\alpha-1}||[A,B_k]u_i||^2.
\end{equation} 

In the same paper, the authors deal with another type of inequalities treated by Harell and Stubbe \cite{HaSt}. For this, we say that a real function $f$ satisfies condition (H1) if there exists a function $r(x)$ such that 
$$\frac{f(x)-f(y)}{x-y}\geq \frac{r(x)+r(y)}{2}\qquad \text{(H1)}.$$
An example of such a function $f$ is whenever $f'$ is concave, in this case $r=f'$. In \cite{AH}, Ashbaugh and Hermi prove the following: 

\begin{Thm} \cite[Thm. 7]{AH}
Under the same assumptions as in the previous theorem, and if $f$ is a function satisfying the condition (H1), we have 
$$\sum_{i=1}^m f(\lambda_i)\left(\sum_{k=1}^N \langle [A,B_k]u_i,B_k u_i\rangle\right)\leq -\frac{1}{2}\sum_{i=1}^m r(\lambda_i)\left(\sum_{k=1}^N ||[A,B_k]u_i||^2\right)+\mathfrak{R}, $$
where 
$$\mathfrak{R}=\sum_{k=1}^N\sum_{i=1}^m\sum_{j=m+1}^\infty |\langle [A,B_k]u_i, u_j\rangle|^2\left(\frac{f(\lambda_i)}{\lambda_{m+1}-\lambda_i}+\frac{1}{2}r(\lambda_i)\right).$$
\end{Thm}

For the particular case, when $f(\lambda)=(\lambda_{m+1}-\lambda)^\alpha$ with $\alpha\geq 2$, they deduce the following inequality \cite[Cor. 8]{AH}
\begin{equation}\label{eq:AH2}
\sum_{i=1}^m\sum_{k=1}^N(\lambda_{m+1}-\lambda_i)^\alpha \langle [A,B_k]u_i,B_k u_i\rangle\leq \frac{\alpha}{2}\sum_{i=1}^m\sum_{k=1}^N (\lambda_{m+1}-\lambda_i)^{\alpha-1}||[A,B_k]u_i||^2.
\end{equation} 

\begin{proof}[Proof of Theorem \ref{thm:recursionformulas}]
In the following we will use Inequalities \eqref{eq:AH1} and \eqref{eq:AH2} in the case of the drifting Hodge Laplacian  defined on a manifold with boundary.  Recall that on a compact Riemannian manifold 
$(\Omega,g)$ with boundary $\partial \Omega$ the Dirichlet problem on differential $p$-forms is given by 
\begin{equation}\label{eq:driftingdirichlet}
\left\{\begin{array}{ll}
\Delta_f\omega=\lambda_{p,f}\omega\,\,  &\text{on}\,\, \Omega\\
\omega=0 \,\,&\text{on}\,\, \partial\Omega.
\end{array}\right.
\end{equation}
We now take $A=\Delta_f$ and $B_k$ a function on $\Omega$, which we will denote by $G$, in Inequalities \eqref{eq:AH1} and \eqref{eq:AH2}. We follow closely the computations done in \cite{IM}. 
Throughout the proof we choose an orthonormal frame \(\{e_i\}_{i=1,\ldots,n}\) of $TM$ such that \(\nabla e_i=0\) at a fixed point.
First, using Equation \eqref{eq:Bochner}, we compute 
\begin{eqnarray*} 
[\Delta_f,G]\omega&=&[\nabla_f^*\nabla,G]\omega\\
&=&\left(-\sum_{i=1}^n\nabla_{e_i}\nabla_{e_i}+\nabla_{df}\right)(G\omega)-G\nabla_f^*\nabla\omega\\
&=&(\Delta_f G)\omega-2\nabla_{dG}\omega.
\end{eqnarray*}
Here, we recall that $\Delta_f G=\Delta G+g(df,dG)$. Therefore, we get that 
\begin{eqnarray} \label{eq:commutatordeltafg}
\int_\Omega\langle [\Delta_f,G]\omega,G\omega\rangle  d\mu_f&=&\int_\Omega G(\Delta_f G)|\omega|^2 d\mu_f
-2\int_\Omega G\langle\nabla_{dG}\omega,\omega\rangle d\mu_f\nonumber\\
&=&\int_\Omega G(\Delta_f G)|\omega|^2 d\mu_f-\frac{1}{2}\int_\Omega \langle dG^2,d(|\omega|^2)\rangle d\mu_f\nonumber\\
&=&\int_\Omega G(\Delta_f G)|\omega|^2 d\mu_f-\frac{1}{2}\int_\Omega |\omega|^2(\Delta_f G^2) d\mu_f\nonumber\\
&=&\int_\Omega |\omega|^2|dG|^2 d\mu_f.
\end{eqnarray}
In the last equality, we use the fact that $\Delta_f G^2=2G\Delta_f G-2|dG|^2$ which can be shown by a straightforward computation. Also, we have that 
\begin{equation}\label{eq:normedeltafg}
\int_\Omega|[\Delta_f,G]\omega|^2 d\mu_f=\int_\Omega (\Delta_f G)^2|\omega|^2 d\mu_f+4\int_\Omega|\nabla_{dG}\omega|^2 d\mu_f-4\int_\Omega\Delta_f G\langle \omega, \nabla_{dG}\omega\rangle d\mu_f.
\end{equation}

In the following, we will assume that the manifold $\Omega$ is a domain in a complete Riemannian manifold $(M^n,g)$ that is isometrically immersed into the Euclidean space $\R^{n+m}$ endowed with its canonical metric.
We choose in the above formulas the function $G$ to be $G=x_l, l=1,\ldots,n+m$ the components of the immersion $X=(x_1,\ldots,x_{n+m})$. 
Recall that $(\partial x_l)^T=d(x_l\circ \iota)$ where $\iota$ is the isometric immersion.  
In addition, $\omega=\omega_i$ are eigenforms of the problem \eqref{eq:driftingdirichlet} associated to the eigenvalues $\lambda_{i,p,f}$ that are chosen to be of $L^2$-norm equal to $1$. The computation done in \eqref{eq:commutatordeltafg} gives after summing over $l$ that 
$$\sum_{l=1}^{m+n}\int_\Omega\langle [\Delta_f,x_l]\omega_i,x_l\omega_i\rangle  d\mu_f
=\sum_{l=1}^{m+n}\int_\Omega |dx_l|^2|\omega_i|^2 d\mu_f=n.$$
Recall here that $\sum_{l=1}^{n+m} |dx_l|^2=n$. By taking $G=x_l$ in Equation \eqref{eq:normedeltafg} and summing over $l$ yields
\begin{eqnarray*}
\sum_{l=1}^{n+m}\int_\Omega|[\Delta_f,x_l]\omega_i|^2 d\mu_f&=&\sum_{l=1}^{n+m}\int_\Omega (\Delta_f x_l)^2|\omega_i|^2 d\mu_f+4\sum_{l=1}^{n+m}\int_\Omega|\nabla_{dx_l}\omega_i|^2 d\mu_f \\&&
-4\sum_{l=1}^{n+m}\int_\Omega(\Delta_f x_l)\langle \omega_i, \nabla_{dx_l}\omega_i\rangle d\mu_f.\\
&=&\sum_{l=1}^{n+m}\int_\Omega \left(\Delta x_l+g(df,dx_l)\right)^2|\omega_i|^2 d\mu_f+4\int_\Omega|\nabla_{dx_l}\omega_i|^2 d\mu_f\\&&
-4\int_\Omega\left(\Delta x_l+g(df,dx_l)\right)\langle \omega_i, \nabla_{dx_l}\omega_i\rangle d\mu_f\\
&=&\int_\Omega (n^2|H|^2+|df|^2+2ng(df,H))|\omega_i|^2 d\mu_f+4\int_\Omega|\nabla\omega_i|^2 d\mu_f\\
&&-2\int_\Omega (n g(H,d|\omega_i|^2)+g(df,d|\omega_i|^2)) d\mu_f\\
&=&\int_\Omega (n^2|H|^2+|df|^2)|\omega_i|^2 d\mu_f+4\int_\Omega|\nabla\omega_i|^2 d\mu_f \\
&&-2\int_\Omega g(df,d|\omega_i|^2) d\mu_f.
\end{eqnarray*} 
Here, we used that the mean curvature of the immersion is given by $H=\frac{1}{n}(\Delta x_1,\ldots, \Delta x_{n+m})$. Applying the Bochner-Weitzenb\"ock formula \eqref{eq:Bochner} to $\omega_i$ and taking the scalar product with 
$\omega_i$ itself, the above equality acquires the form:
\begin{eqnarray*}
\sum_{l=1}^{n+m}\int_\Omega|[\Delta_f,x_l]\omega_i|^2 d\mu_f&=
&\int_\Omega (n^2|H|^2+|df|^2)|\omega_i|^2 d\mu_f \\
&&+4\left(\lambda_{i,p,f}-\frac{1}{2}\int_\Omega\Delta_f(|\omega_i|^2) d\mu_f-\int_\Omega\langle(\mathfrak{B}^{[p]}+T_f^{[p]})\omega_i,\omega_i\rangle d\mu_f\right)\\
&&-2\int_\Omega (\Delta_f f) |\omega_i|^2 d\mu_f\\
&=&\int_\Omega (n^2|H|^2+|df|^2)|\omega_i|^2 d\mu_f \\
&&+4\left(\lambda_{i,p,f}
-\frac{1}{2}\int_{\partial\Omega}\frac{\partial}{\partial\nu}(|\omega_i|^2) d\mu_f
-\int_\Omega\langle(\mathfrak{B}^{[p]}+T_f^{[p]})\omega_i,\omega_i\rangle d\mu_f\right)\\
&&-2\int_\Omega (\Delta f+|df|^2)|\omega_i|^2 d\mu_f\\
&=&\int_\Omega n^2|H|^2|\omega_i|^2 d\mu_f+4\left(\lambda_{i,p,f}-\int_\Omega\langle(\mathfrak{B}^{[p]}+T_f^{[p]})\omega_i,\omega_i\rangle d\mu_f\right)\\&&-\int_\Omega (2\Delta f+|df|^2)|\omega_i|^2 d\mu_f.\\
\end{eqnarray*}
Inserting the above equations into \eqref{eq:AH1} and \eqref{eq:AH2} completes the proof
of Theorem \ref{thm:recursionformulas}.
\end{proof}

\begin{proof}[Proof of Theorem \ref{thm:jacobi}] 
First of all, let $u$ be the function given by 
$$u^{i_1,\ldots,i_{p+1}}=\langle \partial x_{i_1}\wedge\ldots \wedge\partial x_{i_{p+1}},\nu\wedge\omega\rangle$$
where $\omega$ is a $p$-form on $M$ and $\nu$ is the inward unit normal vector field. 
Throughout the proof we will denote by $\nabla^{\R^{n+1}}$ the connection on $\R^{n+1}$ and by $\nabla$ the connection on $M$, as well as for ${d^{\R^{n+1}}}$ and $d$. To simplify the notations, we will denote $u^{i_1,\ldots,i_{p+1}}$ by $u$ in the following technical lemma (see \cite[Lem. 3.1]{IRS} for $p=1$).

\begin{Lem} 
Suppose that $M$ is a $f$-minimal hypersurface of the weighted manifold $(\R^{n+1},g={\rm can}, e^{-f}dv_g)$. Then, we have
\begin{eqnarray*}
L_fu&=&-\langle \partial x_{i_1}\wedge\ldots\wedge \partial x_{i_{p+1}},({\nabla_\nu^{\R^{n+1}}} {d^{\R^{n+1}}}f)^T\wedge\omega\rangle
+2\langle \partial x_{i_1}\wedge\ldots\wedge \partial x_{i_{p+1}},{\bf {\bf II}}(e_i)\wedge\nabla_{e_i}\omega\rangle
\\&&
-\langle \partial x_{i_1}\wedge\ldots\wedge \partial x_{i_{p+1}},\nu\wedge ({\bf {\bf II}}^2)^{[p]}\omega\rangle
+\langle \partial x_{i_1}\wedge\ldots\wedge \partial x_{i_{p+1}},\nu\wedge\Delta_f\omega\rangle\\
&&+\langle \partial x_{i_1}\wedge\ldots\wedge \partial x_{i_{p+1}},\nu\wedge ({\bf {\bf II}}^{[p]})^2\omega\rangle
-\langle \partial x_{i_1}\wedge\ldots\wedge \partial x_{i_{p+1}},\nu\wedge {\widetilde T}_f^{[p]}\omega\rangle \\
&&-{\rm Hess}^{\R^{n+1}} f(\nu,\nu)u.
\end{eqnarray*}
Here ${\widetilde T}_f^{[p]}$ is defined in the same way as in Section 2.2 on $\R^{n+1}$. 
\end{Lem} 

\begin{proof} For any $X\in TM$, we have 
$$X(u)=-\langle \partial x_{i_1}\wedge\ldots \wedge\partial x_{i_{p+1}},{\bf {\bf II}}(X)\wedge\omega\rangle
+\langle \partial x_{i_1}\wedge\ldots \wedge\partial x_{i_{p+1}},\nu\wedge\nabla_X\omega\rangle.$$

Here, we used the fact that $(\nabla^{\R^{n+1}}_X\omega)^T=\nabla_X\omega$, since $\omega$ is a differential form on $M$. Differentiating again with respect to $X$ yields 
\begin{eqnarray*}
X(X(u))&=&-\langle \partial x_{i_1}\wedge\ldots\wedge \partial x_{i_{p+1}},\nabla^{\R^{n+1}}_X {\bf {\bf II}}(X)\wedge\omega\rangle-\langle \partial x_{i_1}\wedge\ldots\wedge \partial x_{i_{p+1}},{\bf {\bf II}}(X)\wedge\nabla^{\R^{n+1}}_X\omega\rangle\\&&-\langle \partial x_{i_1}\wedge\ldots\wedge \partial x_{i_{p+1}}, {\bf {\bf II}}(X)\wedge\nabla_X\omega\rangle+\langle \partial x_{i_1}\wedge\ldots\wedge \partial x_{i_{p+1}}, \nu\wedge\nabla^{\R^{n+1}}_X\nabla_X\omega\rangle\\
&=&-\langle \partial x_{i_1}\wedge\ldots\wedge \partial x_{i_{p+1}},\nabla_X {\bf {\bf II}}(X)\wedge\omega\rangle-|{\bf {\bf II}}(X)|^2\langle \partial x_{i_1}\wedge\ldots\wedge \partial x_{i_{p+1}},\nu\wedge\omega\rangle\\&&-2\langle \partial x_{i_1}\wedge\ldots\wedge \partial x_{i_{p+1}},{\bf {\bf II}}(X)\wedge \nabla_X\omega\rangle-\langle \partial x_{i_1}\wedge\ldots\wedge \partial x_{i_{p+1}}, {\bf {\bf II}}(X)\wedge\nu\wedge {\bf {\bf II}}(X)\lrcorner\omega\rangle\\
&&+\langle \partial x_{i_1}\wedge\ldots\wedge \partial x_{i_{p+1}}, \nu\wedge\nabla_X\nabla_X\omega\rangle.
\end{eqnarray*}
In the last part, we used the Gauss formula and the identity  $\nu\lrcorner\nabla^{\R^{n+1}}_X\omega={\bf {\bf II}}(X)\lrcorner\omega$. Now tracing over an orthonormal frame $\{e_i\}$ of $TM$ gives that 
\begin{eqnarray*} 
\Delta_f u&=&\Delta u+g(df,du)\\
&=&-\langle \partial x_{i_1}\wedge\ldots\wedge \partial x_{i_{p+1}}, (\delta {\bf {\bf II}})\wedge\omega\rangle+|{\bf {\bf II}}|^2\langle \partial x_{i_1}\wedge\ldots\wedge \partial x_{i_{p+1}},\nu\wedge\omega\rangle\\
&&+2\langle \partial x_{i_1}\wedge\ldots\wedge \partial x_{i_{p+1}},{\bf {\bf II}}(e_i)\wedge \nabla_{e_i}\omega\rangle-\langle \partial x_{i_1}\wedge\ldots\wedge \partial x_{i_{p+1}},\nu\wedge ({\bf {\bf II}}^2)^{[p]}\omega\rangle\\
&&-\langle \partial x_{i_1}\wedge\ldots\wedge \partial x_{i_{p+1}}, \nu\wedge\nabla_{e_i}\nabla_{e_i}\omega\rangle-\langle \partial x_{i_1}\wedge\ldots\wedge \partial x_{i_{p+1}},{\bf {\bf II}}(df)\wedge\omega\rangle\\&&+\langle \partial x_{i_1}\wedge\ldots\wedge \partial x_{i_{p+1}},\nu\wedge\nabla_{df}\omega\rangle.\\
&=&n\langle \partial x_{i_1}\wedge\ldots\wedge \partial x_{i_{p+1}}, dH\wedge\omega\rangle+|{\bf {\bf II}}|^2u+2\langle \partial x_{i_1}\wedge\ldots\wedge \partial x_{i_{p+1}},{\bf {\bf II}}(e_i)\wedge \nabla_{e_i}\omega\rangle\\&&-\langle \partial x_{i_1}\wedge\ldots\wedge \partial x_{i_{p+1}},\nu\wedge ({\bf {\bf II}}^2)^{[p]}\omega\rangle+\langle \partial x_{i_1}\wedge\ldots\wedge \partial x_{i_{p+1}}, \nu\wedge\nabla_f^*\nabla\omega\rangle\\&&-\langle \partial x_{i_1}\wedge\ldots\wedge \partial x_{i_{p+1}},{\bf {\bf II}}(df)\wedge\omega\rangle.
\end{eqnarray*}
In the last equality, we used the fact that $\delta {\bf {\bf II}}=-ndH$ and the expression of $\nabla_f^*\nabla$. Now, combining Equation \eqref{eq:bp} with the Bochner-Weitzenb\"ock formula \eqref{eq:Bochner} combined with \eqref{eq:bp} that is shown in the next section, the above equality reduces to
\begin{eqnarray*} 
\Delta_f u&=&n\langle \partial x_{i_1}\wedge\ldots\wedge \partial x_{i_{p+1}}, dH\wedge\omega\rangle+|{\bf {\bf II}}|^2u+2\langle \partial x_{i_1}\wedge\ldots\wedge \partial x_{i_{p+1}},{\bf {\bf II}}(e_i)\wedge \nabla_{e_i}\omega\rangle\\&&-\langle \partial x_{i_1}\wedge\ldots\wedge \partial x_{i_{p+1}},\nu\wedge ({\bf {\bf II}}^2)^{[p]}\omega\rangle+\langle \partial x_{i_1}\wedge\ldots\wedge \partial x_{i_{p+1}}, \nu\wedge\Delta_f\omega\rangle\\&&
+\langle \partial x_{i_1}\wedge\ldots\wedge \partial x_{i_{p+1}}, \nu\wedge \left(({\bf {\bf II}}^{[p]})^2-nH {\bf {\bf II}}^{[p]}-T_f^{[p]}\right)\omega\rangle \\
&&-\langle \partial x_{i_1}\wedge\ldots\wedge \partial x_{i_{p+1}},{\bf {\bf II}}(df)\wedge\omega\rangle.
\end{eqnarray*}
Now, for a $f$-minimal hypersurface, we write 
$$d^{\R^{n+1}}f=df+\frac{\partial f}{\partial\nu}\nu=df-nH\nu.$$
Hence, by differentiating along a vector $X\in TM$, we get that $$(\nabla^{\R^{n+1}}_X d^{\R^{n+1}}f)^T=\nabla_X df+nH {\bf {\bf II}}(X).$$ 
Therefore, we find
$$(\widetilde T_f^{[p]}\omega)^T=\sum_{i=1}^n e_i\wedge (\nabla^{\R^{n+1}}_{e_i} d^{\R^{n+1}}f)^T\lrcorner\omega=T_f^{[p]}\omega+nH {\bf {\bf II}}^{[p]}\omega.$$

Also, we have that 
$$ndH=-d\left(\frac{\partial f}{\partial\nu}\right)=-(\nabla^{\R^{n+1}}_\nu d^{\R^{n+1}}f)^T+{\bf {\bf II}}(df).$$ 
Hence, plugging these computations into the expression of $\Delta_f u$, we arrive at
\begin{eqnarray*}
\Delta_f u&=&-\langle \partial x_{i_1}\wedge\ldots\wedge \partial x_{i_{p+1}}, (\nabla^{\R^{n+1}}_\nu 
{d^{\R^{n+1}}}f)^T \wedge\omega\rangle \\
&&+|{\bf {\bf II}}|^2u+2\langle \partial x_{i_1}\wedge\ldots\wedge \partial x_{i_{p+1}},{\bf {\bf II}}(e_i)\wedge \nabla_{e_i}\omega\rangle\\&&-\langle \partial x_{i_1}\wedge\ldots\wedge \partial x_{i_{p+1}},\nu\wedge ({\bf {\bf II}}^2)^{[p]}\omega\rangle+\langle \partial x_{i_1}\wedge\ldots\wedge \partial x_{i_{p+1}}, \nu\wedge\Delta_f\omega\rangle\\&&+\langle \partial x_{i_1}\wedge\ldots\wedge \partial x_{i_{p+1}}, \nu\wedge ({\bf {\bf II}}^{[p]})^2\omega\rangle-\langle \partial x_{i_1}\wedge\ldots\wedge \partial x_{i_{p+1}}, \nu\wedge {\widetilde T}_f^{[p]}\omega\rangle.
\end{eqnarray*}
The proof of the lemma then follows from the expression of $L_f$. 
\end{proof}

In order to complete the proof of Theorem \ref{thm:jacobi}, we proceed as in \cite{IRS}. Let $\{\varphi_j\}$ be an $L^2(e^{-f}dv_g)$-orthonormal basis of eigenfunctions of $L_f$ associated to $\lambda_j(L_f)$. Let $E^d$ be the sum of the first $d$-eigenspaces of $\Delta_f$: 
$$E^d=\mathop\oplus_{j=1}^d V_{\Delta_f}(\lambda_{j,p,f}).$$

For all $l\geq 2$, we consider the following system of equations 
$$\int_M u^{i_1,\ldots,i_{p+1}}\varphi_1 d\mu_f=\ldots=\int_M u^{i_1,\ldots,i_{p+1}}\varphi_{l-1} d\mu_f=0,$$
where we recall that $u^{i_1,\ldots,i_{p+1}}=\langle \partial x_{i_1}\wedge\ldots \partial x_{i_{p+1}}, \nu\wedge \omega\rangle$. This system consists of $\begin{pmatrix}n+1\\p+1\end{pmatrix}(l-1)$ homogeneous linear equations in the variable $\omega\in E^d$. Hence, if we take $d= \begin{pmatrix}n+1\\p+1\end{pmatrix}(l-1)+1$, we can find a non-trivial $\omega\in E^d$ which is orthogonal to the first $(l-1)$-eigenfunctions of $L_f$. 
Therefore, we deduce that 
$$\lambda_l(L_f)\int_M (u^{i_1,\ldots,i_{p+1}})^2 d\mu_f\leq 
\int_M u^{i_1,\ldots,i_{p+1}} L_f(u^{i_1,\ldots,i_{p+1}})d\mu_f.$$

Summing over $i_1<\ldots<i_{p+1}$, we first have for the l.h.s. that 
$$\sum_{i_1<\ldots<i_{p+1}} (u^{i_1,\ldots,i_{p+1}})^2=|\nu\wedge\omega|^2=|\omega|^2.$$ 
For the r.h.s. of the previous inequality, we use the previous lemma to get that 
\begin{eqnarray*}
\sum_{i_1<\ldots<i_{p+1}} u^{i_1,\ldots,i_{p+1}} L_f(u^{i_1,\ldots,i_{p+1}})&=&2\langle \nu\wedge\omega,{\bf {\bf II}}(e_i)\wedge\nabla_{e_i}\omega\rangle-\langle  \nu\wedge\omega,\nu\wedge ({\bf {\bf II}}^2)^{[p]}\omega\rangle\\&&+\langle  \nu\wedge\omega,\nu\wedge\Delta_f\omega\rangle+\langle  \nu\wedge\omega,\nu\wedge ({\bf {\bf II}}^{[p]})^2\omega\rangle\\&&-\langle  \nu\wedge\omega,\nu\wedge {\widetilde T}_f^{[p]}\omega\rangle
-{\rm Hess}^{\R^{n+1}} f(\nu,\nu)|\nu\wedge\omega|^2\\
&=&-\langle  ({\bf {\bf II}}^2)^{[p]}\omega,\omega\rangle+\langle  \Delta_f\omega,\omega\rangle+\langle ({\bf {\bf II}}^{[p]})^2\omega,\omega\rangle\\&&-\langle  {\widetilde T}_f^{[p]}\omega,\omega\rangle-
{\rm Hess}^{\R^{n+1}} f(\nu,\nu)|\omega|^2.
\end{eqnarray*}
Hence, we deduce that 
\begin{eqnarray*}
\lambda_l(L_f)\int_M |\omega|^2 d\mu_f&\leq& \int_M \big(-\langle  ({\bf {\bf II}}^2)^{[p]}\omega,\omega\rangle+\langle  \Delta_f\omega,\omega\rangle+\langle ({\bf {\bf II}}^{[p]})^2\omega,\omega\rangle-\langle  {\widetilde T}_f^{[p]}\omega,\omega\rangle \\
&&-{\rm Hess}^{\R^{n+1}} f(\nu,\nu)|\omega|^2\big)d\mu_f\\
&\leq & \int_M \left(-\langle  ({\bf {\bf II}}^2)^{[p]}\omega,\omega\rangle+\langle  \Delta_f\omega,\omega\rangle+\langle ({\bf {\bf II}}^{[p]})^2\omega,\omega\rangle-(p+1)a|\omega|^2\right)d\mu_f.
\end{eqnarray*}
In the last equality, we used the fact that any eigenvalue of the operator ${\widetilde T}_f^{[p]}$ is the sum of $p$-distinct eigenvalues of ${\rm Hess}^{\R^{n+1}}f$. Now, we have that 
$$\int_M\langle  \Delta_f\omega,\omega\rangle d\mu_f\leq \lambda_{d(l),p,f}\int_M|\omega|^2d\mu_f,$$
since $\omega\in E^{d(l)}$. Also, we have that 
\begin{eqnarray*}
-\langle  ({\bf {\bf II}}^2)^{[p]}\omega,\omega\rangle+\langle ({\bf {\bf II}}^{[p]})^2\omega,\omega\rangle&=&-\sum_{i=1}^n\langle {\bf {\bf II}}(e_i)\wedge {\bf II}(e_i)\lrcorner\omega,\omega\rangle+|\sum_{i=1}^n e_i\wedge {\bf II}(e_i)\lrcorner\omega|^2\\
&=&-\sum_{i,j=1}^n\langle e_j\lrcorner {\bf II}(e_i)\lrcorner\omega,e_i\lrcorner {\bf II}(e_j)\lrcorner\omega\rangle\\
&\leq & \gamma_M p(p-1)|\omega|^2.
\end{eqnarray*} 
Finally, we deduce that 
$$\lambda_l(L_f) \leq \lambda_{d(l),p,f}-(p+1)a+\gamma_M p(p-1),$$
which finishes the proof of the theorem.
\end{proof}

\begin{proof}[Proof of Corollaries \ref{cor:jacobi} and \ref{cor:jacobibis}]
Let $l_0$ be the largest integer such that $d(l_0)\leq \beta$. Thus, 
$${\rm Ind}_f(M)\geq l_0.$$
Now as we have that $d\left(\frac{1}{\begin{pmatrix}n+1\\p+1\end{pmatrix}}\beta+1-\frac{1}{\begin{pmatrix}n+1\\p+1\end{pmatrix}}\right)\leq \beta$, then  $$l_0\geq\left[\frac{1}{\begin{pmatrix}n+1\\p+1\end{pmatrix}}\beta+1-\frac{1}{\begin{pmatrix}n+1\\p+1\end{pmatrix}}\right]\geq \frac{1}{\begin{pmatrix}n+1\\p+1\end{pmatrix}}\beta.$$ This proves the first corollary. To prove the second one, we take as before $$l_0=\left[\frac{1}{\begin{pmatrix}n+1\\p+1\end{pmatrix}}b_p(M)+1-\frac{1}{\begin{pmatrix}n+1\\p+1\end{pmatrix}}\right],$$ then, we clearly have that $d(l_0)\leq b_p(M)$. Therefore, $\lambda_{d(l_0),p,f}=0$. In the case of a self-shrinker, we recall that for $f=\frac{|X|^2}{2}$ the hypersurface $M$ is $f$-minimal and that 
${\rm Hess}^{\R^{n+1}} f=g={\rm can}$. Hence with the assumption $p(p-1)\gamma_M\leq p+1$, we deduce that $\lambda_{l_0}(L_f)\leq 0$. Finally, we get that the index is at least $n+1+l_0$ by the fact that $L_f$ has at least $n+1$ eigenvalues equal to $-1$. This finishes the proof.
\end{proof}

\section{Geometric applications of the main results}
In this section we will provide several geometric applications
of Theorem \ref{thm:recursionformulas} by making explicit choices of the function \(f\).\\

First, we will consider the case  when $f$ is the Riemannian distance function in order to compute explicitly the different terms in the statement of Theorem \ref{thm:recursionformulas}.
To this end, let $M^n\to \R^{n+m}$ be an isometric immersion and $\Omega$ a domain in $M$. Consider the function $f:\Omega\to \mathbb{R}$ given by $f(X)=a\frac{|X|^2}{2}$, where $a$ is a real positive number. This particular choice of \(f\) has important applications in mean curvature flow
as described in the introduction of the article. 
Here, $|\cdot|$ denotes the Euclidean norm in $\R^{n+m}$ and $X=(x_1,\ldots,x_{n+m})$ are the components of the immersion. In other words, the function $f$ is the square of the distance function from the origin point $0\in \R^{n+m}$ to a point $X\in \Omega$.  Using the decomposition $X=X^T+X^\perp$, we have that $d^{\R^{n+m}}f=a X$ and therefore, $df=a X^T$. Hence, we deduce that
\begin{equation}\label{eq:df}
|df|^2=a^2(|X|^2-|X^\perp|^2).
\end{equation}
For any $Y\in TM$, we compute
\begin{eqnarray} \label{eq:nablahypersurface}
\nabla_Ydf&=&a \nabla_YX^T \nonumber \\
&=&a (\nabla^{\R^{n+m}}_YX^T)^T\nonumber\\
&=&a(Y-\nabla^{\R^{n+m}}_YX^\perp)^T\nonumber\\
&=&a\left(Y-\sum_{i=1}^n(\nabla^{\R^{n+m}}_YX^\perp,e_i)e_i\right)\nonumber\\
&=&a\left(Y+\sum_{i=1}^n(X^\perp,\nabla^{\R^{n+m}}_Ye_i)e_i\right)\nonumber\\
&=&a\left(Y+\sum_{i=1}^n(X^\perp,{\bf II}(Y,e_i))e_i\right)\nonumber\\
&=&a(Y+{\bf II}_{X^\perp}(Y)),
\end{eqnarray}
where \({\bf II}_{X^\perp}\) is defined as in the proof of Theorem \ref{thm:inequality}.
By tracing \eqref{eq:nablahypersurface}, we deduce that 
\begin{equation}\label{eq:Deltaf}
\Delta f=a(-n-n(X,H)).
\end{equation}
Also, plugging \eqref{eq:nablahypersurface} into the expression of $T_f^{[p]}$ yields
$$T_f^{[p]}\omega=a(p\omega+{\bf II}^{[p]}_{X^\perp}\omega),$$
where ${\bf II}^{[p]}_{Z}$ is the canonical extension of ${\bf II}_{Z}$ to $p$-forms 
as defined in \eqref{dfn:sff-forms}.
Hence, we deduce 
\begin{equation}\label{eq:tpf}
\langle T_f^{[p]}\omega,\omega\rangle=a\left(p|\omega|^2+\langle {\bf II}^{[p]}_{X^\perp}\omega,\omega\rangle\right).
\end{equation}
In the following, we will bound the term $\langle\mathfrak{B}^{[p]}\omega_i,\omega_i\rangle $ in Theorem \ref{thm:recursionformulas} by using the results of \cite{S3}. For this, recall that A. Savo shows in \cite[Thm. 1]{S3} that for any isometric immersion $(M^n,g)\to (N^{n+m},g)$, the Bochner operator $\mathfrak{B}^{[p]}$ on $p$-forms of $M$ splits as 
$$\mathfrak{B}^{[p]}=\mathfrak{B}^{[p]}_{\rm ext}+\mathfrak{B}^{[p]}_{\rm res},$$
where $\mathfrak{B}^{[p]}_{\rm ext}$ is the operator defined by 
$$\mathfrak{B}^{[p]}_{\rm ext}=\sum_{j=1}^{m}\left({\rm trace}({\bf II}_{\nu_j}){\bf II}_{\nu_j}^{[p]}-{\bf II}_{\nu_j}^{[p]}\circ {\bf II}_{\nu_j}^{[p]}\right).$$ 
Here, $\{\nu_1,\ldots,\nu_{m}\}$ is a local orthonormal frame of $TM^\perp$, and the operator  $\mathfrak{B}^{[p]}_{\rm res}$ is the operator that satisfies 
$$p(n-p)\gamma_N\leq\mathfrak{B}^{[p]}_{\rm res}\leq p(n-p)\Gamma_N,$$
where $\gamma_N$ and $\Gamma_N$ are respectively a lower bound and an upper bound for the curvature operator of $N$. Hence for an isometric immersion $M\to (\R^{n+m},{\rm can})$, we deduce that 
\begin{equation}\label{eq:bp}
\mathfrak{B}^{[p]}=\mathfrak{B}^{[p]}_{\rm ext}={\bf II}_{nH}^{[p]}-\sum_{j=1}^{m}{\bf II}_{\nu_j}^{[p]}\circ {\bf II}_{\nu_j}^{[p]}.
\end{equation}
Inserting Equations \eqref{eq:df}, \eqref{eq:Deltaf}, \eqref{eq:tpf} and \eqref{eq:bp} in Theorem \ref{thm:recursionformulas}, we deduce 

\begin{Cor} \label{eq:distancefunction0} Let $X:(M^n,g)\to (\mathbb{R}^{n+m},{\rm can})$ be an isometric immersion and let $f=a\frac{|X|^2}{2}$ where $a$ is a real positive number. For any $p\in\{0,\ldots,n\}$, the eigenvalues of the drifting Hodge Laplacian $\Delta_f$ acting on $p$-forms on a domain $\Omega$ of $M$ with Dirichlet boundary conditions satisfy for any $k\geq 1$
\begin{eqnarray*}
\sum_{i=1}^k(\lambda_{k+1,p,f}-\lambda_{i,p,f})^\alpha&\leq& \frac{4}{n}\sum_{i=1}^k(\lambda_{k+1,p,f}-\lambda_{i,p,f})^{\alpha-1}\big(\lambda_{i,p,f}
-\int_\Omega\langle {\bf II}_{(aX^\perp+nH)}^{[p]}\omega_i,\omega_i\rangle d\mu_f \\
&&+\sum_{j=1}^{m}\int_\Omega |{\bf II}_{\nu_j}^{[p]}\omega_i|^2d\mu_f-ap+\frac{na}{2} \\
&&-\frac{a^2}{4}\int_\Omega |X|^2|\omega_i|^2d\mu_f+\frac{1}{4}\int_\Omega |aX^\perp+nH|^2|\omega_i|^2d\mu_f\big),
\end{eqnarray*}
for $\alpha\leq 2$. 
Also, we have 
\begin{eqnarray*}
\sum_{i=1}^k(\lambda_{k+1,p,f}-\lambda_{i,p,f})^\alpha&\leq& \frac{2\alpha}{n}\sum_{i=1}^k(\lambda_{k+1,p,f}-\lambda_{i,p,f})^{\alpha-1}\big(\lambda_{i,p,f}
-\int_\Omega\langle {\bf II}_{(aX^\perp+nH)}^{[p]}\omega_i,\omega_i\rangle d\mu_f  \\
&&+\sum_{j=1}^{m}\int_\Omega |{\bf II}_{\nu_j}^{[p]}\omega_i|^2d\mu_f-ap+\frac{na}{2}\\
&&-\frac{a^2}{4}\int_\Omega |X|^2|\omega_i|^2d\mu_f+\frac{1}{4}\int_\Omega |aX^\perp+nH|^2|\omega_i|^2d\mu_f\big),
\end{eqnarray*}
for $\alpha\geq 2$. 
\end{Cor}

The result in Corollary \ref{eq:distancefunction0} generalizes the one in \cite[Thm. 1.1]{Z} when taking $a=\frac{1}{2}, \alpha=2$ and $\Omega=M$ being a domain in $\R^m
$. In this case, ${\bf II}=0$ and we get that  
\begin{eqnarray*}
\sum_{i=1}^k(\lambda_{k+1,p,f}-\lambda_{i,p,f})^2&\leq& \frac{4}{n}\sum_{i=1}^k(\lambda_{k+1,p,f}-\lambda_{i,p,f})\left(\lambda_{i,p,f}-\frac{p}{2}+\frac{n}{4}-\frac{1}{16}\mathop{\rm min}\limits_\Omega (|X|^2)\right).
\end{eqnarray*}

Moreover, Corollary \ref{eq:distancefunction0} generalizes the result in \cite[Thm 1.2]{Z} when taking for some integer $l>1$, the number $a=l-1, \alpha=2$ and $(M^n,g)=(\R^{n-l}\times \mathbb{S}^l(1),\langle,\rangle_{\R^{n-l}}\oplus\langle,\rangle_{\mathbb{S}^l})$ where $\langle,\rangle_{\mathbb{S}^l}$ is the standard metric on the unit round sphere $\mathbb{S}^l(1)$   of curvature $1$. Indeed, we take the immersion $M\hookrightarrow \R^{n-l}\times \R^{l+1}$ where the second fundamental form is given by the matrix
$\begin{pmatrix}
0&0\\
0&{\rm Id}
\end{pmatrix}$. Thus, we have  $X^\perp=-\nu$ and $nH={\rm trace}({\bf II})=l\nu$ and, therefore, $aX^\perp+nH=\nu$. Also, using that ${\bf II}^{[p]}\omega=\sum_{i=1}^n e_i\wedge {\bf II}(e_i)\lrcorner \omega$ for any $p$-form $\omega$, we deduce by a straightforward computation that 
$$-\langle {\bf II}^{[p]}\omega,\omega\rangle+|{\bf II}^{[p]}\omega|^2=\sum_{i,j=n-l+1}^n |e_i\lrcorner e_j\lrcorner \omega|^2\leq p(p-1)|\omega|^2.$$ 
Hence, we get that
\begin{align*}
\sum_{i=1}^k(\lambda_{k+1,p,f}-\lambda_{i,p,f})^2\leq & \frac{4}{n}\sum_{i=1}^k(\lambda_{k+1,p,f}-\lambda_{i,p,f})\big(\lambda_{i,p,f}+\frac{2n(l-1)+1-4p(l-1)+4p(p-1)}{4}\\&-\frac{(l-1)^2}{4}\mathop{\rm min}\limits_\Omega (|X|^2)\big).
\end{align*} 

Corollary \ref{eq:distancefunction0} also generalizes the one in \cite{CP} on compact self-shrinkers. In this case, for $a=1$, we get that 

\begin{eqnarray*}
\sum_{i=1}^k(\lambda_{k+1,p,f}-\lambda_{i,p,f})^2&\leq& \frac{4}{n}\sum_{i=1}^k(\lambda_{k+1,p,f}-\lambda_{i,p,f})\big(\lambda_{i,p,f}+\sum_{j=1}^{m}\int_\Omega |{\bf II}_{\nu_j}^{[p]}\omega_i|^2d\mu_f\\&&-p+\frac{n}{2}-\frac{1}{4}\mathop{\rm min}\limits_\Omega (|X|^2)\big).
\end{eqnarray*}

In the last part of this section, we will consider the case when $f$ is the distance function from some fixed point in $M$ with respect to the Riemannian metric $g$ on $M$. For this, fix a point $x_0\in M$ and consider the distance function 
$$d_{x_0}:\Omega\to [0,\infty[, d_{x_0}(x)=d(x_0,x)$$
and $\rho_{x_0}(x):=\frac{1}{2}d_{x_0}(x)^2$. We recall from \cite{Pe} that the distance function $d_{x_0}$ is smooth in the complement of the cut locus of $x_0$ and that its gradient is of norm $1$ almost everywhere. Thus we will choose a domain $\Omega\subset M$ such that it is contained in this complement. Let us recall the comparison theorem \cite{Pe} (see also \cite{HaSi}). For this, we consider for any $l\in \mathbb{R}$, the function 
$$
H_l(r)=\left\{\begin{array}{ll}
(n-1)\sqrt{l}\,{\rm cot}(\sqrt{l}r),\,\,\,\,  l>0\\\\
\frac{n-1}{r}, \,\,\,\,l=0\\\\
(n-1)\sqrt{|l|}\,{\rm coth}(\sqrt{|l|}r),\,\,\,\, l<0.
\end{array}\right.
$$

\begin{Thm}  \cite{Pe} \label{thm:comparison} Let $(M^n,g)$  be a complete Riemannian manifold. 
\begin{enumerate}
\item If ${\rm Ric}^M\geq (n-1)l$ for some $l\in \mathbb{R}$, then for every $x_0\in M$, the inequalities
$$\Delta d_{x_0}(x)\geq -H_l(d_{x_0}(x)), \quad\text{and}\quad \Delta \rho_{x_0}(x)\geq -(1+d_{x_0}(x)H_l(d_{x_0}(x)),$$
hold at smooth points of $d_{x_0}$. Moreover, these inequalities hold on the whole manifold in the sense of distributions. 
\item If the sectional curvature satisfies $l_1\leq K^M\leq l_2$ and $\gamma$ is a minimizing geodesic starting from $x_0\in M$ such that its image is disjoint to the cut locus of $x_0$, then 
$$\nabla^2 d_{x_0}(X,X)\leq \frac{H_{l_1}(t)}{n-1}g(X,X), \quad\text{and}\quad \nabla^2 d_{x_0}(X,X)\geq \frac{H_{l_2}(t)}{n-1}g(X,X),$$
for $X\perp \dot{\gamma}(t)$ and 
$\nabla^2d_{x_0}(\dot{\gamma}(t),\dot{\gamma}(t))=0,$
for $t\in [0,L]$ for some $L$. As a consequence, one has that 
$$\nabla^2 \rho_{x_0}(X,X)\leq \frac{tH_{l_1}(t)}{n-1}g(X,X), \quad\text{and}\quad \nabla^2 \rho_{x_0}(X,X)\geq \frac{tH_{l_2}(t)}{n-1}g(X,X),$$
for $X\perp \dot{\gamma}(t)$ and 
$\nabla^2\rho_{x_0}(\dot{\gamma}(t),\dot{\gamma}(t))=1$.
\end{enumerate}
\end{Thm}

We will now use Theorem \ref{thm:comparison} to compute the different terms in Theorem \ref{thm:recursionformulas}. We begin with the case when ${\rm Ric}^M\geq (n-1)l$ for some $l\in \R$ on the manifold $M$. By taking $f=a \rho_{x_0}$ where $a>0$, we deduce that the inequalities
$$|df|^2=a^2 d_{x_0}^2,\quad\text{and}\quad \Delta f(x)\geq -a\left(1+d_{x_0}(x)H_l(d_{x_0}(x)\right)$$
hold in the sense of distributions. Therefore, we deduce from Theorem \ref{thm:recursionformulas} the following

\begin{Cor} \label{thm:recursionformulas1} Let $X:(M^n,g)\to (\mathbb{R}^{n+m},{\rm can})$ be an isometric immersion. Assume that ${\rm Ric}^M\geq (n-1)l$ for some $l\in \mathbb{R}$. Let $x_0$ be  fixed point in $M$ and $\Omega$ a domain in the complement of the cut locus of $x_0$. Let $f=a\rho_{x_0}$ for some positive real number $a>0$. The eigenvalues of the drifting Hodge Laplacian $\Delta_f$ acting on $p$-forms on $\Omega$ with Dirichlet boundary conditions satisfy for any $k\geq 1$
\begin{eqnarray*}
\sum_{i=1}^k(\lambda_{k+1,p,f}-\lambda_{i,p,f})^\alpha&\leq& \frac{4}{n}\sum_{i=1}^k(\lambda_{k+1,p,f}-\lambda_{i,p,f})^{\alpha-1}\big(\lambda_{i,p,f}-\delta_1+\frac{1}{4}\delta_2\big),
\end{eqnarray*}
for $\alpha\leq 2$, where $\delta_1=\mathop{\rm inf}\limits_\Omega (\mathfrak{B}^{[p]}+T_f^{[p]})$ and $\delta_2$ is given by
$$\delta_2=\mathop {\rm sup}\limits_{x\in\Omega}\left(n^2|H|^2+2a(1+d_{x_0}(x)H_l(d_{x_0}(x))-a^2d_{x_0}(x)^2\right).$$ 
Also, we have 
\begin{eqnarray*}
\sum_{i=1}^k(\lambda_{k+1,p,f}-\lambda_{i,p,f})^\alpha&\leq& \frac{2\alpha}{n}\sum_{i=1}^k(\lambda_{k+1,p,f}-\lambda_{i,p,f})^{\alpha-1}\big(\lambda_{i,p,f}-\delta_1+\frac{1}{4}\delta_2\big),
\end{eqnarray*}
for $\alpha\geq 2$. 
\end{Cor}

In the following, we will consider the case when the sectional curvature of $M$ is bounded from below by $l_1$ and from above by $l_2$. We let $\delta_2$ be the number given by 
\small
$$\delta_2:=\left\{\begin{array}{ll}
\mathop{\rm sup}\limits_{x\in \Omega}\left(n^2|H|^2-4a\left((p-1)\frac{d_{x_0}(x)H_{l_2}(d_{x_0}(x))}{n-1}+1\right)+2a(1+d_{x_0}H_{l_1}(d_{x_0}))-a^2 d_{x_0}^2\right),\,\,\,\,  l_2\leq 0\\\\
\mathop{\rm sup}\limits_{x\in \Omega}\left(n^2|H|^2-4ap\frac{d_{x_0}(x)H_{l_2}(d_{x_0}(x))}{n-1}+2a(1+d_{x_0}H_{l_1}(d_{x_0}))-a^2 d_{x_0}^2\right), \,\,\,\,l_2>0.\\\\
\end{array}\right.
$$
\normalsize
We have

\begin{Cor} Let $X:(M^n,g)\to (\mathbb{R}^{n+m},{\rm can})$ be an isometric immersion. Assume that $l_1\leq K^M\leq l_2$ for some $l_1,l_2\in \mathbb{R}$. Let  $\Omega$ a domain in $M$ such that $\Omega$ is contained in the complement of the cut locus of $x_0\in \Omega$. Let $f=a\rho_{x_0}$ for some positive real number $a>0$. For any $p\in\{0,\ldots,n\}$, the eigenvalues of the drifting Hodge Laplacian $\Delta_f$ acting on $p$-forms on a domain $\Omega$ of $M$ with Dirichlet boundary conditions satisfy for any $k\geq 1$
\begin{eqnarray*}
\sum_{i=1}^k(\lambda_{k+1,p,f}-\lambda_{i,p,f})^\alpha&\leq& \frac{4}{n}\sum_{i=1}^k(\lambda_{k+1,p,f}-\lambda_{i,p,f})^{\alpha-1}\big(\lambda_{i,p,f}-\int_\Omega\langle {\bf II}_{nH}^{[p]}\omega_i,\omega_i\rangle d\mu_f \\&&
+\sum_{j=1}^{m-n}\int_\Omega |{\bf II}_{\nu_j}^{[p]}\omega_i|^2d\mu_f+\frac{1}{4}\delta_2
\end{eqnarray*}
for $\alpha\leq 2$. Also, we have 
\begin{eqnarray*}
\sum_{i=1}^k(\lambda_{k+1,p,f}-\lambda_{i,p,f})^\alpha&\leq& \frac{2\alpha}{n}\sum_{i=1}^k(\lambda_{k+1,p,f}-\lambda_{i,p,f})^{\alpha-1}\big(\lambda_{i,p,f}-\int_\Omega\langle {\bf II}_{nH}^{[p]}\omega_i,\omega_i\rangle d\mu_f\\&&
+\sum_{j=1}^{m-n}\int_\Omega |{\bf II}_{\nu_j}^{[p]}\omega_i|^2d\mu_f+\frac{1}{4}\delta_2
\end{eqnarray*}
for $\alpha\geq 2$. 
\end{Cor}

\begin{proof} We proceed as in \cite{S1}. At any point $x\in \Omega\setminus\{x_0\}$, there is an orthonormal frame $\{e_1(x),\ldots, e_{n-1}(x),e_n=\nabla d_{x_0}(x)\}$ such that 
$\nabla^2\rho_{x_0}$ has the eigenvalues $\eta_1(x),\ldots,\eta_{n-1}(x),1$. By the comparison theorem \ref{thm:comparison}, we get that for any $j=1,\ldots,n-1$
$$\frac{d_{x_0}(x)H_{l_2}(d_{x_0}(x))}{n-1}\leq \eta_j(x)\leq \frac{d_{x_0}(x)H_{l_1}(d_{x_0}(x))}{n-1}.$$
Therefore, we find that 
$$|df|^2=a^2 d_{x_0}^2 \quad\text{and}\quad\Delta f\geq -a(1+d_{x_0}H_{l_1}(d_{x_0})).$$
Recall that the endomorphism $T^{[p]}_f$ is by definition the sum of $p$ distinct eigenvalues of $\nabla^2 f$. Hence, for any $p$-form $\omega$, we get that 
$$\langle T^{[p]}_f\omega,\omega\rangle\geq \left\{\begin{array}{ll}
a\left((p-1)\frac{d_{x_0}(x)H_{l_2}(d_{x_0}(x))}{n-1}+1\right)|\omega|^2,\,\,\,\,  l_2\leq 0\\\\
ap\frac{d_{x_0}(x)H_{l_2}(d_{x_0}(x))}{n-1}|\omega|^2, \,\,\,\,l_2>0.\\
\end{array}\right.
$$
This is because of $\frac{d_{x_0}(x)H_{l_2}(d_{x_0}(x))}{n-1}\geq 1$ if $l_2\leq 0$ and $\frac{d_{x_0}(x)H_{l_2}(d_{x_0}(x))}{n-1}\leq 1$ if $l_2>0$. Replacing all these inequalities in Theorem \ref{thm:recursionformulas}, we get the result after combining with Equation \eqref{eq:bp} for the curvature term. 
\end{proof}


\begin{thebibliography}{50} 
\bibitem[Asa]{Asa} S. Asada,
{\it On the first eigenvalue of the {L}aplacian acting on {$p$}-forms},
Hokkaido Math. J., \textbf{9(1)} (1980), 112-122.

\bibitem[Ash]{Ash} M. Ashbaugh, \textit{Isoperimetric and universal inequalities for eigenvalues.}, Spectral theory and geometry (Edinburgh, 1998), London Math. Soc. Lecture Note Ser., \textbf{273}, Cambridge Univ. Press, Cambridge, 1999. 

\bibitem[AH]{AH} M. S. Ashbaugh, L. Hermi,  {\it On Harell-Stubbe type inequalities for the discrete spectrum of a self-adjoint operator}, arxiv:0712.4396

\bibitem[BE]{BE} D. Bakry, M. Emery, {\it Diffusions hypercontractives}, S\'eminaire de probabilit\'es, XIX, 1983/84,
Lecture Notes in Math., \textbf{1123}, Springer, Berlin, 1985. 

\bibitem[BO]{BO} J. Baldauf, T. Ozuch,
{\it Spinors and mass on weighted manifolds},
Commun. Math. Phys. \textbf{394}, 1153–1172 (2022). 

\bibitem[BK]{BK} J. Baldauf, D. Kim,
{\it Parabolic frequency on Ricci flows},
Int. Math. Res. Not. IMRN 2023, no. \textbf{12}, 10098-10114.

\bibitem[BCP]{BCP} M. Batista, M.P. Cavalcante and J. Pyo, {\it Some isoperimetric inequalities and eigenvalue estimates in weighted manifolds}, J. Math. Anal. Appl., \textbf{419} (2014), 617-626.

\bibitem[CGH]{CGH} F. Chami, N. Ginoux and G. Habib, {\it New eigenvalue estimates involving Bessel functions}, Publ. Mat. \textbf{65} (2021), 681-726.

\bibitem[CC]{CC} D. Chen, Q.-M. Cheng, {\it Extrinsic estimates for eigenvalues of the Laplace operator}, 
J. Math. Soc. Japan \textbf{60 (2)} (2008) 325--339

\bibitem[CP]{CP} Q.-M. Cheng, Y. Peng,  {\it Estimates for eigenvalues of $\mathfrak{L}$ operators on self-shrinkers}, Comm. Contemp. Math., \textbf{15} (2013), 135011.

\bibitem[CY]{CY} Q.-M. Cheng, H. Yang, {\it Bounds on eigenvalues of Dirichlet Laplacian},
Math. Ann. \textbf{337} (2007), no. 1, 159–175

\bibitem[CZ]{CZ} D. Chen and Y. Zhang, {\it Inequalities for eigenvalues of the weighted
Hodge Laplacian}, arXiv:1312.0218.

\bibitem[CM]{CM} T. Colding, W. Minicozzi,
{\it Generic mean curvature flow I: generic singularities},
Ann. of Math. (2) \textbf{175} (2012), no. 2, 755–833. 

\bibitem[EHI]{EHI} A. El Soufi, E. Harell and S. Ilias
{\it Universal inequalities for the eigenvalues
of Laplace and Schr\"odinger operators on
submanifolds},
Trans. AMS \textbf{361} (2009), 2337–2350. 

\bibitem[GM]{GM} S. Gallot, D. Meyer, {\it Op\'erateur de courbure et Laplacien des formes diff\'erentielles d'une vari\'et\'e riemannienne}, J. Math. Pures Appl. \textbf{54} (1975), 259-284. 

\bibitem[GS]{GS} P. Gu\'erini, A. Savo, {\it Eigenvalue and gap estimates for the Laplacian acting on $p$-forms}, Trans. Amer. Math. Soc. \textbf{356} (2003), 319-344.

\bibitem[HaSt]{HaSt} E. M. Harell, J. Stubbe, {\it On trace identities and universal eigenvalue estimates for some partial differential operators}, Trans. Amer. Math. Soc. \textbf{349} (1997), 1797-1809. 

\bibitem[HaSi]{HaSi} A. Hassanezhad, A. Siffert, {\it A note on Kuttler-Sigillito's inequalities},  Ann. Math. Qué. \textbf{44}(1), 125–147 (2020)

\bibitem[HP]{HP} G.N. Hile, M.H. Protter, {\it Inequalities for eigenvalues of the Laplacian},  Indiana Univ. Math. J. \textbf{29} (1980) 523--538.

\bibitem[IM]{IM} S. Ilias, O. Makhoul, {\it Universal inequalities for the eigenvalues of the Hodge de Rham Laplacian}, Ann. Glob. Anal. Geom. \textbf{36} (2009), 191-204.

\bibitem[IRS]{IRS} D. Impera, M. Rimoldi and A. Savo, {\it Index and first Betti number of $f$-minimal hypersurfaces and self-shrinkers}, Rev. Mat. Ibero. \textbf{36} (2019), 817-840.

\bibitem[JMZ]{JMZ}  M. Carmen Domingo-Juan, V. Miquel and J. Zhu, {\it Reilly's type inequality for the Laplacian associated to a density related with shrinkers for MCF}, J. Diff. Equat. \textbf{272} (2021), 
958-978. 

\bibitem[LM]{LM} H.B. Lawson, M.-L. Michelsohn, {\it Spin geometry},
Princeton Mathematical Series, 38. Princeton University Press, Princeton, NJ, 1989. 

\bibitem[Li]{Li} A. Lichnerowicz, {\it Spineurs harmoniques}, C. R. Acad. Sci. Paris, \textbf{257} (1963), 7-9.

\bibitem[L]{L} J. Lott, {\it Some geometric properties of the Bakry-Emery Ricci tensor}, Comment. Math. Helv. \textbf{78} (2003), 865-883.

\bibitem[MD]{MD} L. Ma, S.H. Du, {\it Extension of Reilly formula with applications to eigenvalue estimates for drifting Laplacians}, C. R. Math. Acad. Sci. Paris \textbf{348} (2010), 1203-1206.


\bibitem[ML]{ML} L. Ma, B.Y. Liu, {\it Convexity of the first eigenfunction of the drifting Laplacian operator and its applications}, New York J. Math. \textbf{14} (2008), 393-401.


\bibitem[ML1]{ML1} L. Ma, B.Y. Liu, {\it Convex eigenfunction of a drifting Laplacian operator and the fundamental gap}, Pacific J. Math. \textbf{240} (2009), 343-361.


\bibitem[PGW1]{PGW1} L.E. Payne, G. Polya, H.F. Weinberger, {\it Sur le quotient de deux fréquences propres cons\'ecutives}, 
C. R. Acad. Sci., \textbf{Paris 241} (1955) 917--919.

\bibitem[PGW2]{PGW2} L.E. Payne, G. Polya, H.F. Weinberger, 
{\it On the ratio of consecutive eigenvalues}, J. Math. Phys. \textbf{35 (1)} (1956) 289--298. 


\bibitem[Pe]{Pe} P. Petersen, {\it Riemannian geometry}, Graduate Texts in Mathematics. Springer, New York, second edition, 2006.

\bibitem[PW]{PW} P. Petersen, M. Wink, {\it The Bochner technique and weighted curvatures}, SIGMA \textbf{16} (2020), 064.

\bibitem[Per]{Per} G.Perelman, {\it The entropy formula for the Ricci flow and its geometric applications}, arXiv:math/0211159

\bibitem[P]{P} D. Le Peutrec, {\it On Witten Laplacians and Brascamp-Lieb's inequality on manifolds with boundary}, Integ. Equa. Opera. Theo. \textbf{87} (2017), 411-434. 


\bibitem[S]{S} A. Savo, {\it On the first Hodge eigenvalue of isometric immersions}, Proc. Amer. Math. Soc. \textbf{133} (2004), 587-594.

\bibitem[S1]{S1} A. Savo, {\it On the lowest eigenvalue of the Hodge Laplacian on compact negatively curved domains}, Ann. Glob. Anal. Geom. \textbf{35} (2009), 39-62. 

\bibitem[S2]{S2} A. Savo, {\it Index bounds for minimal hypersurfaces of the sphere}, Indiana Univ. Math. J. \textbf{36} (2010), 823-837. 

\bibitem[S3]{S3} A. Savo, {\it The Bochner formula for isommetric immersions}, Pac. J. Math.  \textbf{272} (2014), 395-422.

\bibitem[XX]{XX} C. Xia, H. Xu, {\it Inequalities for eigenvalues of the drifting Laplacian on Riemannian manifolds}, Ann. Glob. Anal. Geom. \textbf{45} (2014), 155-166.

\bibitem[Y]{Y} H.C. Yang, {\it An estimate of the difference between consecutive eigenvalues}, Preprint \textbf{IC/91/60} of ICTP, Trieste, 1991.

\bibitem[Z]{Z1} L. Zeng, {\it Eigenvalues of the drifting Laplacian on complete non-compact Riemannian manifolds}, Nonlinear Anal. \textbf{141} (2016), 1-15.

\bibitem[Z1]{Z} L. Zeng, {\it Estimates for the eigenvalues of the drifting Laplacian on some complete Ricci solitons}, Kyushu J. Math. \textbf{72} (2018), 143-156. 



\end{thebibliography}
\end{document}